\documentclass{article}

\PassOptionsToPackage{hyphens}{url}\usepackage[colorlinks=red,linkcolor=green,citecolor=blue]{hyperref}
\usepackage{amsmath,amssymb,amsfonts,amsthm}
\usepackage[authoryear,round]{natbib}
\usepackage{graphicx}
\usepackage{color}
\usepackage{mathtools}
\usepackage{subfig}
\usepackage{multirow}
\usepackage{float}
\usepackage{appendix}
\usepackage{booktabs}
\usepackage{longtable}
\usepackage{cleveref}
\Crefname{equation}{}{}
\usepackage{tikz}
\usetikzlibrary{decorations.pathreplacing,arrows}
\usepackage[margin=1in]{geometry}
\usepackage{algorithmic}
\usepackage{algorithm}

\allowdisplaybreaks
\newcommand{\mc}[1]{\mathcal{#1}}
\crefname{appsec}{Appendix}{Appendices}

\title{A Node Formulation for Multistage Stochastic Programs with Endogenous Uncertainty}

\author{Giovanni Pantuso\\
  Department of Mathematical Sciences\\
  University of Copenhagen, 2100 Copenhagen, Denmark\\
  \href{mailto:gp@math.ku.dk}{gp@math.ku.dk}
}
\date{\today}

\begin{document}
\maketitle

\begin{abstract}
  This paper introduces a node formulation for multistage stochastic programs with
  endogenous (i.e., decision-dependent) uncertainty. Problems with such structure
  arise when the choices of the decision maker determine a change in the likelihood of future
  random events. The node formulation avoids an explicit statement of non-anticipativity constraints,
  and as such keeps the dimension of the model sizeable.
  An exact solution algorithm for a special case is introduced and tested on a case study.
  Results show that the algorithm outperforms a commercial solver as the size of the instances increases.
\end{abstract}

\section{Introduction}\label{sec:introduction}

Multistage stochastic programs offer a viable framework for modeling and solving problems involving a sequence
of decisions interspersed with partial resolutions of some stochastic process. At each decision stage the
decision maker knows the content of the uncertainty resolved until that stage and a probabilistic characterization
of the remaining stochastic process. In the classical settings, it is assumed that decisions do not modify the
stochastic process in any way \citep[see e.g.,][]{KalW94,BirL97}. In other words, the uncertainty is entirely exogenous.
While this description fits a large number of decision problems, several other can be found where decisions have
an influence on the remainder stochastic process, e.g., by changing the likelihood of future realizations.
This category of problem is referred to as \textit{multistage stochastic programs with endogenous uncertainty} (MSPEU).

Following \citep{GoeG06}, there exist at least two ways in which decisions can influence the underlying stochastic process.
The first possibility is that decisions alter the probability distribution of the stochastic process,
thus changing the likelihood of the possible events. The second possibility is that decisions determine the
time when the uncertainty is (partially) resolved. This article is concerned with multistage stochastic programs
affected by the first type of endogenous uncertainty.

The research dealing with decisions influencing probability distributions is rather sparse. \citet{JonWW98}
consider a case where decisions influence both the probability measure and the timing of the observation,
i.e., at which stage certain random variables will be observed. The framework includes both two-stage and
multistage problems, though the decisions influencing the uncertainty must be made at the first-stage only.  
\citet{Ahm00} illustrates examples of problems with endogenous uncertainty, such as facility location, network design and
server selection. The author presents an exact solution method for the resulting one-stage integer problems.
\citet{VisSF04} consider the problem of investing in strengthening actions for the links of a network subject to
disruptive events. The problem is modeled as a two-stage stochastic program where first-stage investment decisions
influence the likelihood of disruptive events happening at upgraded links. The same problem is studied also by
\citet{Fla10}, \citet{PeeSGV10} and \citet{LauPK14}. \citet{HelW05} consider the problem of interdicting a
stochastic network, that is a network whose structure is unknown to the interdictor. In this problem, the
probabilities of different future network configurations depend on previous interdiction actions.
\citet{TonFR12} present an oil refinery planning problem considering that the uncertainty in product
yield is influenced by operation mode changeovers. \citet{Hel16} and \citet{HelBT18} discuss several ways of incorporating
the influence of decision variables on the underlying probability distributions in two-stage stochastic
programs. Particularly, the authors formulate two-stage models where prior probabilities are distorted
through an affine transformation, or combined using a convex combination of several probability distributions.
Furthermore, the authors present models which incorporate the  parameters of the probability distribution as
first-stage decision variables. Finally, \citet{EscGMU18} study the problem of mitigating the effects of
natural disasters through preventive actions. The problem is formulated as a three-stage stochastic
bilinear integer program with both exogenous and endogenous uncertainty. Particularly, decisions can influence
both the probabilities and the intensity of future uncertain events. 

When the second type of uncertainty is considered, the decision making process is such that the uncertainty
is not resolved automatically at each decision stage as in classical stochastic programs with only exogenous
uncertainty. Rather, the decisions made implicitly determine the time when the uncertainty is resolved. A typical
example, provided by \citet{GoeG04} and based on the petrochemical industry is as follows. A decision maker has
to decide which gas reservoirs to explore, and when, by installing exploration facilities. The size and quality
of the reservoirs is uncertain and can be known only when facilities have been installed. Thus, the time the
uncertainty is resolved depends on installation decisions. The literature dealing with this type of stochastic
programs includes \citet{ColM08}, \citet{TarG08}, \citet{TarGG09}, \citet{ColM10}, \citet{GupG11}, \citet{MerV11},
\citet{TarGG13}, \citet{ApaG16}.

The contribution of this paper is as follows. First, we introduce a novel node-formulation for multistage stochastic
programs where decisions (at all stages) influence probability distributions at future stages. A node formulation automatically ensures
non-anticipative decisions and thus avoids writing explicit non-anticipativity constraints (NACs) that are necessary
when the uncertainty is represented via scenarios. NACs have typically been addressed as a bottleneck of available models
for MSPEU, and motivated recent research to find and discard redundant NACs, see \citet{ApaG16}, \citet{HooM16} and \cite{MooM18}.
The new formulation is based on a novel scenario tree structure which incorporates the possibility that several (finitely many)
different distributions for later stages can emanate based on the decisions made at a certain decision stage. The new scenario tree structure represents
the second contribution of this paper. Third, we propose a general-purpose efficient algorithm for a special class of MSPEUs.
We demonstrate the use of the new formulation, scenario tree structure and solution algorithm on instances of the \textit{Football Team Composition Problem}.
The instances are made available online for the benefit of future research at \url{https://github.com/GioPan/instancesFTCPwithEndogenousUncertainty}.

The remainder of this paper is organized as follows. In \Cref{sec:MSPEU} we provide a model formulation and
new scenario tree structure for MSPEUs. In \Cref{sec:algorithm} we describe a solution algorithm for a special case of
MSPEUs. In \Cref{sec:cs} we present a case study in which we formulate the Football Team Composition Problem as a MSPEU
and solve it using our specialized algorithm. Finally, we provide concluding remarks in \Cref{sec:conclusions}.

\section{Multistage stochastic programs with decision-dependent uncertainty}\label{sec:MSPEU}
The decision maker is concerned with a sequence of decisions $(x_{t})_{t=1}^T$ at decision stages $t=1,\ldots,T$,
conditional on the partial resolution of a random process $(\xi_t)_{t=1}^T$. At decision stage $t$, decisions are non-anticipative,
meaning that they are based only on the realization of the random process up to, and including, $\xi_t$. The realization of the
remaining random process $\xi_{t+1},\ldots,\xi_T$ is still uncertain. For $t=2,\ldots,T$, the probability
distribution of $\xi_{t},\ldots,\xi_T$ is dependent on past decisions $x_1,\ldots,x_{t-1}$. The resulting
multistage stochastic program is thus characterized by endogenously defined uncertainty.

\begin{figure}
\includegraphics[width=0.7\textwidth]{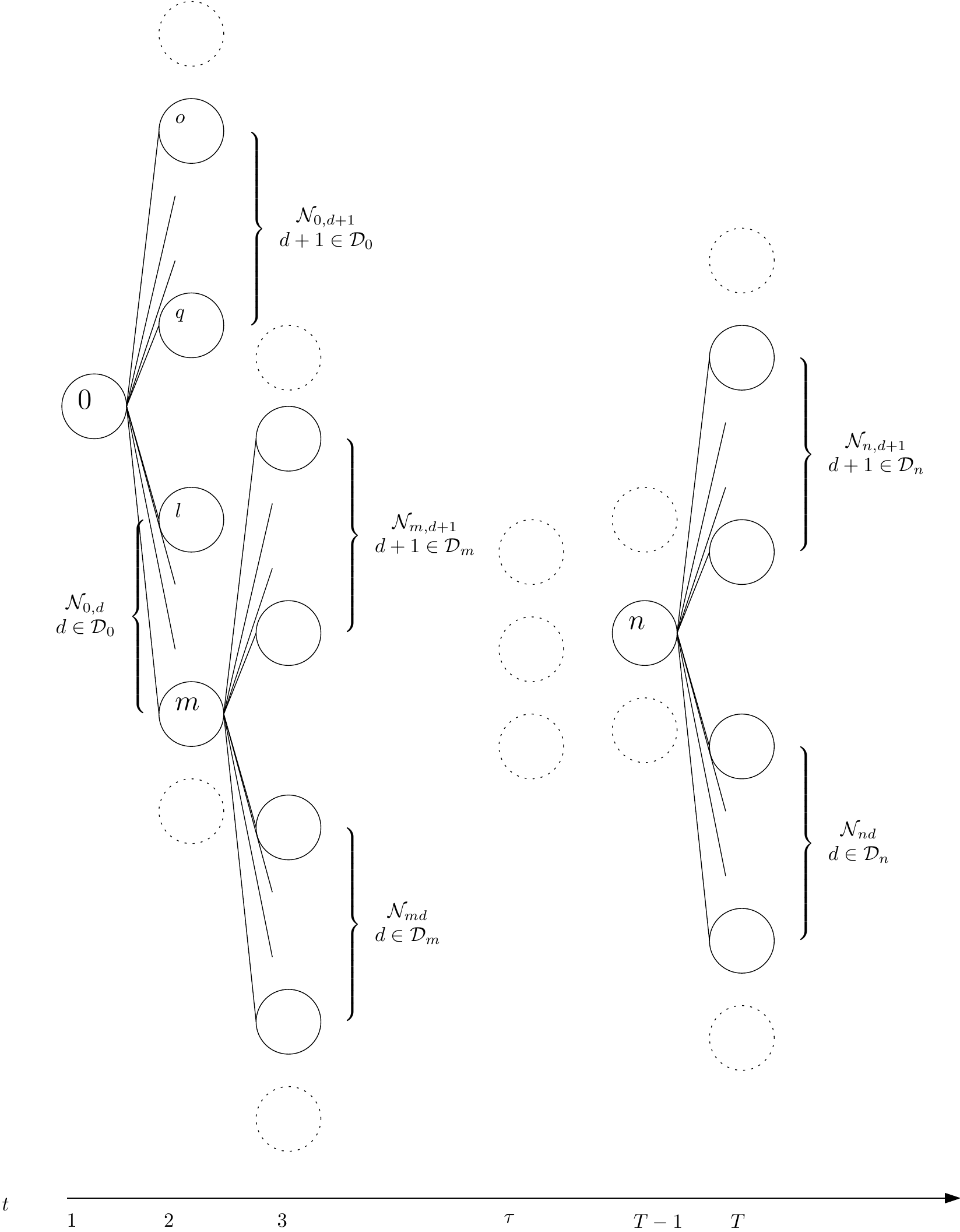}
\caption{Multidistribution scenario tree.}\label{fig:multiDistributionTree}
\end{figure}

Consistently with \citet{JonWW98} we assume that the set of potential probability distributions enforced by
decisions is finite and countable. Furthermore, we assume that the probability distributions are discrete
(possibly after a scenario generation phase). The latter assumption is rather standard and required in order
to solve real-life stochastic programs (except perhaps for a number of specific applications).
The resulting decision-dependent discrete stochastic process can be depicted by means of the scenario tree
structure illustrated in \Cref{fig:multiDistributionTree}. In what follows we refer
to this scenario tree structure as a \textit{multi-distribution scenario tree} (MDST) (see, e.g., also how \citet{KauMWTHF14} modified the classical structure of scenario trees in order to account for multiple time resolutions).

In an MDST the root node, arbitrarily named $0$, represents the current state of the world, when first-stage decisions are made.
First-stage decisions will enforce one out of a finite number of distributions represented by the set $\mc{D}_0$. As an example, in \Cref{fig:multiDistributionTree},
decisions might determine, among other, distributions $d\in\mc{D}_0$ or $d+1\in\mc{D}_0$.  Distribution $d$ is characterized by realizations, represented by nodes,
that include $l$ and $m$, while distribution $d+1$ is characterized by realizations that include $q$ and $o$. 
Similarly, at stage $t=2$, assuming realization $m$ of distribution $d$ occurs at stage $t=1$, the actions of the decision maker will determine
one out of a number of different distributions including $d$ and $d+1$ from the set $\mc{D}_m$. This process continues in a similar manner
until stage $T-1$. Given a realization $n$ at stage $T-1$ the actions of the decision maker will determine one out of a finite number of distributions
$\mc{D}_n$, which include $d$ and $d+1$. Finally, at stage $T$ all the uncertainty is resolved and the decision maker makes final decisions. 

\begin{figure}
\includegraphics[width=0.6\textwidth]{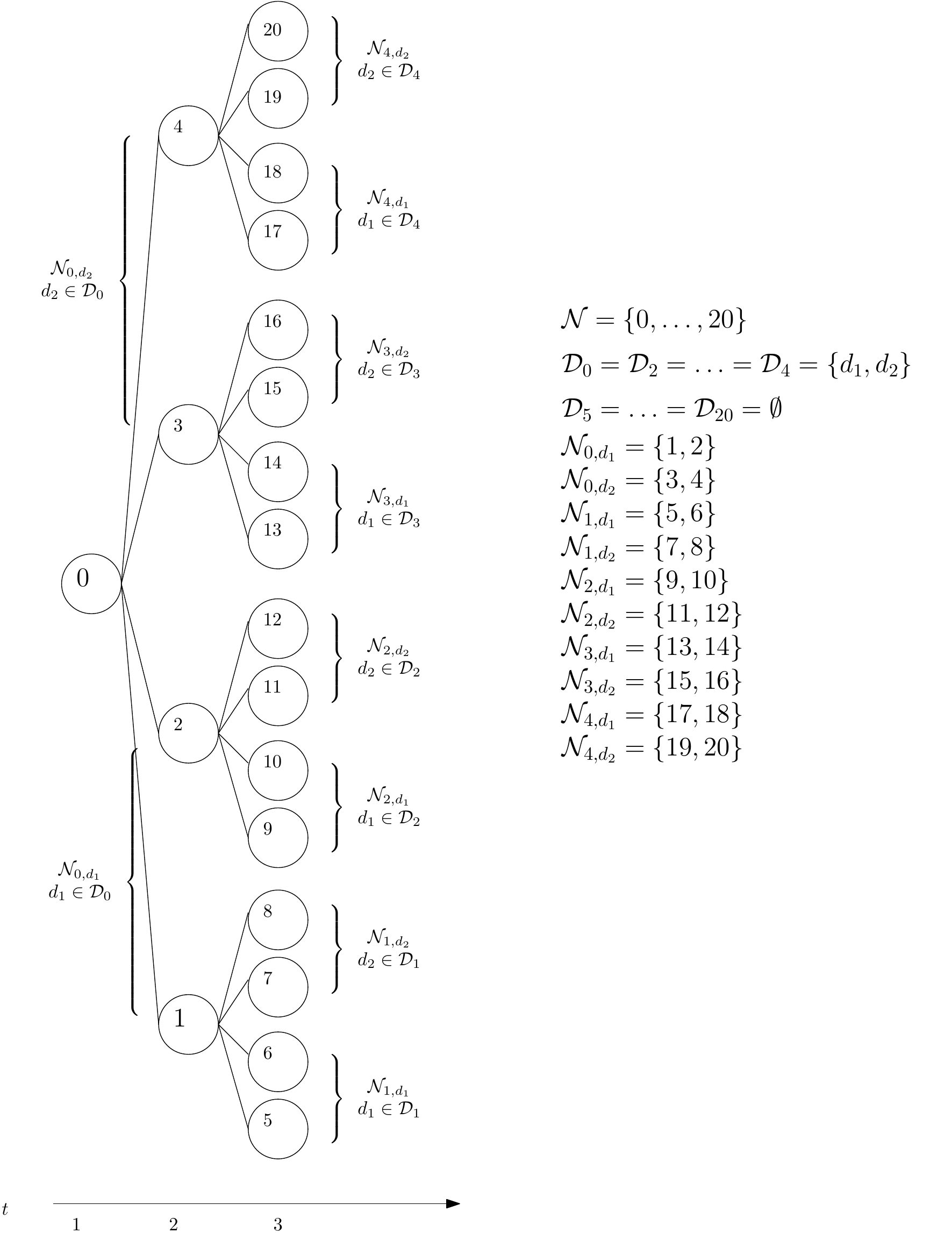}
\caption{Example multidistribution scenario tree with three stages, two possible distributions for each node, and two possible realizations for each distribution.}\label{fig:multiDistributionTreeExample}
\end{figure}

Given an MDST, let us introduce the notation necessary for formulating the MSPEU.
Let $\mathcal{N}$ be the set of nodes in the scenario tree, $\mathcal{N}_t$ be the set of nodes at stage $t$,
$0$ the root node, $t(n)\in\{1,\ldots,T\}$ the stage of node $n$, and $a(n)\in\mathcal{N}$ the parent node
of node $n$ except the root node. Let $\mathcal{D}_n$ be the set of possible distributions which can be
enforced by decision made at node $n\in\mathcal{N}$, and $\mathcal{N}_{nd}$ be the child nodes of node
$n$ if distribution $d\in \mathcal{D}_n$ is enforced. Let $\pi_n$ be the probability of node $n$ with $\pi_0=1$,
and $\sum_{m\in\mathcal{N}_{dn}}\pi_m = \pi_n$ for all $n\in \mathcal{N}$, $d\in \mathcal{D}_n$.
An example of this notation is provided in \Cref{fig:multiDistributionTreeExample} for a three-stage scenario tree,
with two possible distributions emanating from each node and each distribution being characterized by two possible realizations.

Let decision variables $x_n\in \mathbb{R}^{N_{t(n)}}$, $n\in\mathcal{N}$ represent decisions made at node $n$.
These decisions may be integer or fractional and represent ordinary decision made in the decision process.
Let $\delta_{nd}$ be a binary variable which captures the probability distribution enforced by the decisions made at node $n$.
It takes value $1$ if the decisions at node $n$ enforce probability distribution $d$ at the
child nodes of $n$, $0$ otherwise. Finally, let $\theta_n\in \mathbb{R}^1$ be a decision variable which holds the expected value of the decisions made at the nodes descending from $n$.
For the reader's convenience, the notation is also reported in \Cref{app:notation} in a tabular format.  An MSPEU is thus 

\begin{subequations}\label{eq:eumsp}
  \begin{align}
    \label{eq:eumsp:obj} \max~& r_0^Tx_{0}+\sum_{\mathclap{d\in \mathcal{D}_0}}q_{0d}\delta_{0d}+\theta_0  \\
    \label{eq:eumsp:c1}  \text{s.t.}&\sum_{\mathclap{d\in \mathcal{D}_n}}\delta_{nd} = 1 &n\in \mathcal{N},\\
    \label{eq:eumsp:c2}   &A_nx_n + \sum_{\mathclap{d\in \mathcal{D}_n}}B_{nd}\delta_{nd}+ C_{a(n)}x_{a(n)}   + \sum_{\mathclap{d\in \mathcal{D}_{a(n)}}}D_{a(n),d}\delta_{a(n),d}= h_n& n\in \mathcal{N},\\
    \label{eq:eumsp:c3}   &\theta_n = \sum_{d\in \mathcal{D}_n}\delta_{nd}\bigg(\sum_{m\in \mathcal{N}_{nd}}\pi_m(r_m^Tx_{m}+\sum_{\mathclap{d\in\mathcal{D}_m}}q_{md}\delta_{md}+\theta_m)\bigg) & n\in \mathcal{N}\setminus{\mathcal{N}_T},\\
    \label{eq:eumsp:c4}    &\theta_{n} = \Theta_n & n\in \mathcal{N}_{T},\\
    \label{eq:eumsp:c5} & x_n \in X_{t(n)} & n\in \mathcal{N},\\
    \label{eq:eumsp:c6} & \delta_{nd} \in \{0,1\} & n\in \mathcal{N},d\in \mathcal{D}_n,\\
    \label{eq:eumsp:c7} & \theta_{n} \in \mathbb{R} & n\in \mathcal{N}.
  \end{align}
\end{subequations}

where, for each $n\in\mc{N}$, $r_n\in \mathbb{R}^{N_{t(n)}}$ and $q_{nd}\in \mathbb{R}^1$ represent the rewards of decisions $x_n$ and $\delta_{nd}$, respectively,
$A_n\in \mathbb{R}^{M_{t(n)}\times N_{t(n)}}$, $B_{nd}\in \mathbb{R}^{M_{t(n)}\times 1}$, $C_{a(n)}\in \mathbb{R}^{M_{t(n)}\times N_{t(a(n))}}$,
and $D_{a(n),d}\in \mathbb{R}^{M_{t(n)}\times 1}$ are matrices of coefficients and $h_n\in \mathbb{R}^{M_{t(n)}}$ a right-hand-side vector,
with the assumption that $C_{a(0),d}\coloneqq\mathbf{0}$ and $D_{a(0),d}\coloneqq\mathbf{0}$.
Finally, $\Theta_n \in \mathbb{R}^1$ represents the future expected value at leaf node $n$, and $X_{t(n)}\subseteq \mathbb{R}^{N_{t(n)}}$
the domain of the $x_n$ variables.  Objective function \eqref{eq:eumsp:obj}
represents the sum of the profit for the decisions made a the root node ($n=0$) and expected profit of future decisions.
 Constraints \eqref{eq:eumsp:c1} ensure that the decisions made at each node determine exactly one probability distribution.
That is, the stochastic phenomena following the decisions at node $n$, are described by exactly one probability distribution (among the available ones),
enforced by the decisions made.
Constraints \cref{eq:eumsp:c2} describe the dependency between decision variables $x$ and $\delta$ at
each node, and between them and the corresponding decision variables at the parent node.
That is, the choice of a probability distribution at node $n$, $\delta_{nd}$, depends on the decisions $x_n$ made at the same node as well as on
the decisions $x_{a(n)}$ made at the parent node and on the consequent choice of a probability distribution, $\delta_{a(n),d}$. 
Constraints \cref{eq:eumsp:c3} ensure that decision variables $\theta_n$ hold the expected value of future decisions calculated according to the probability chosen.
Consider a generic node $n$ other than a leaf node, and note that, according to \cref{eq:eumsp:c1}, at this node there will be exactly one $\delta_{nd}$ equal to one, that is exactly one distribution be chosen. Thus, the right-hand-side of constraints \cref{eq:eumsp:c3} will be equal to the term of the outer summation corresponding to the index $d$ whose $\delta_{nd}$ is set to one.
The remaining terms are zero.
Correspondingly, $\theta_n$ will hold the expected value calculated according to the chosen probability distribution $d$. 
Constraints \cref{eq:eumsp:c3} can be linearized using standard techniques. A linear reformulation is provided in \Cref{sec:app:bigMref}
and a general procedure to determine the necessary big-$M$ values is described in \Cref{sec:app:bigMs}.
Constraints \cref{eq:eumsp:c4} set the future expectation at the leaf nodes. Finally, constraints \cref{eq:eumsp:c5,eq:eumsp:c6,eq:eumsp:c7} set the
domain of the decision variables. Particularly, $X_{t(n)}$ represents the domain of the $x_n$ variables and may impose integrality restrictions on some/all variables.

\begin{figure}
\includegraphics[width=0.4\textwidth]{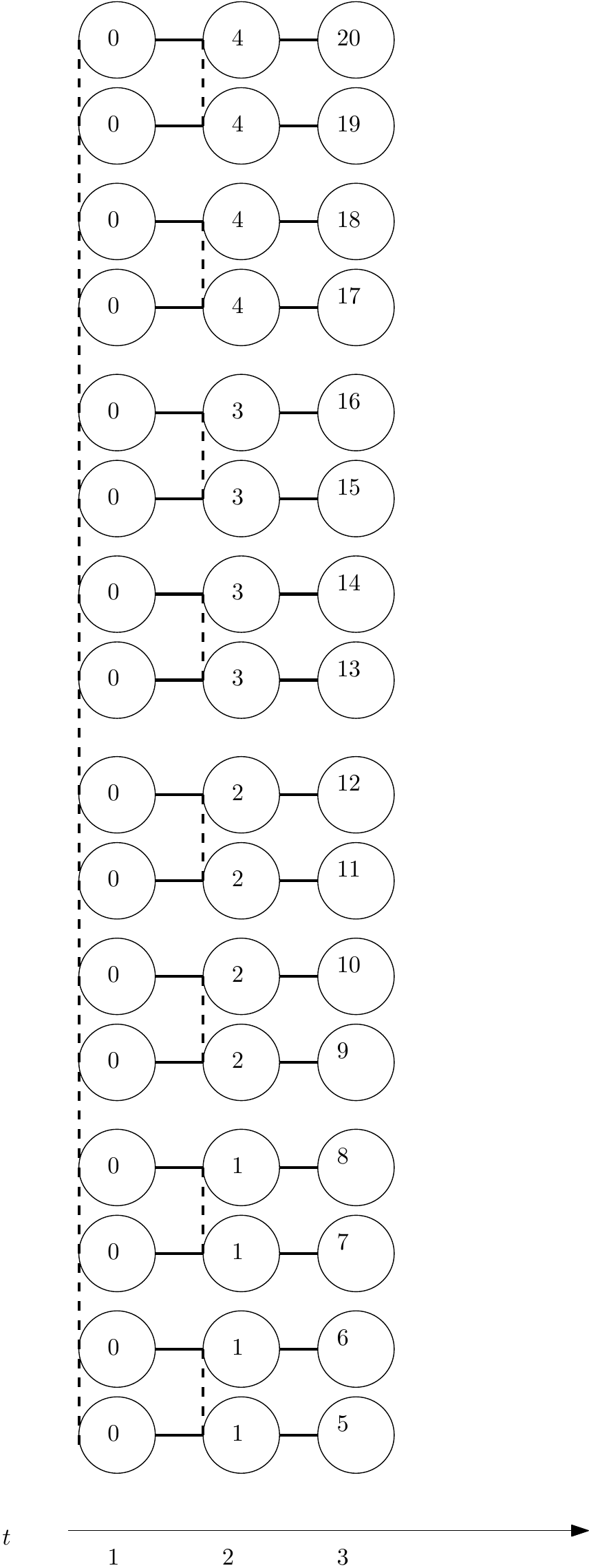}
\caption{Scenarios in the example multidistribution scenario tree in \Cref{fig:multiDistributionTreeExample}. Plain lines connect the nodes belonging to the same scenario. Dashed lines represent non-anticipativity constraints.}\label{fig:multiDistributionTreeExampleScen}
\end{figure}

A node formulation implicitly includes non-anticipativity, that is, automatically ensures that the decisions made at a given stage
are only based on available information. On the other hand, a classical scenario formulation requires non-anticipativity explicitly enforced by means
of constraints, and in general, generates a much larger problem. This can be illustrated by the following  example. Consider
\Cref{fig:multiDistributionTreeExampleScen} which provides the scenario representation of the example MDST in \Cref{fig:multiDistributionTreeExample}.
The number of scenarios is 15, and is the same as the number of leaf nodes in \Cref{fig:multiDistributionTreeExample}.
Assuming the decision at each stage are represented by $N$ decision variables and have to satisfy $M$ constraints, the corresponding scenario formulation would include
\begin{itemize}
\item $N\times 15\times 3$ decision variables, where $3$ is the number of stages,
\item $M\times 15\times 3$ constraints, and
\item approximately $N\times (15+8)$ non-anticipativity constraints (approximately $15$ for the first stage and $8$ for the second stage, though more efficient specifications may be possible).
\end{itemize}
The node formulation would include
\begin{itemize}
\item $N\times 21$ variables (where $20$ is the number of nodes in the scenario tree in \Cref{fig:multiDistributionTreeExample}), and
\item $M\times 21$ constraints.
\end{itemize}
Clearly, node formulations generate, in general, much smaller problems with non-anticipativity constraints playing an important role in scaling up the dimension of a scenario formulation, see e.g.,
\citet{ApaG16}, \citet{HooM16} and \cite{MooM18} for how to reduce the number of non-anticipativity constraints.

\section{A solution algorithm for a special case}
\label{sec:algorithm}

In this section we introduce an algorithm for solving the special case of problem \eqref{eq:eumsp} with $C_{a(n)} = \mathbf{0}$ for all $n$.
The algorithm builds on the fact that, when $C_{a(n)}= \mathbf{0}$, the link between decisions at subsequent stages is created by the $\delta$ variables,
i.e., those that enforce a probability distribution for the next stage. In this case, since we have finitely many distributions, we are allowed to enumerate the expected values obtainable at each node in the scenario tree. This procedure would become impractical when $C_{a(n)}$ is different from $\mathbf{0}$, unless further assumptions on $x_n$ are made.

Let us thus consider the following equivalent formulation of problem \eqref{eq:eumsp}.
\begin{subequations}\label{eq:eumsp_ref}
  \begin{align}
    \label{eq:eumsp_ref:obj} z=\max~& r_0^Tx_{0}+\sum_{d\in \mathcal{D}_0}q_{0d}\delta_{0d}+\theta_0  \\
    \label{eq:eumsp_ref:c1}  \text{s.t.}&\sum_{d\in \mathcal{D}_n}\delta_{nd} = 1 &n\in \mathcal{N},\\
    \label{eq:eumsp_ref:c2}   &A_nx_n + \sum_{d\in \mathcal{D}_n}B_{nd}\delta_{nd} + \sum_{d\in \mathcal{D}_{a(n)}}D_{a(n),d}\delta_{a(n),d}= h_n& n\in \mathcal{N},\\
    \label{eq:eumsp_ref:c3}   &\theta_n\leq \phi_{nd} + M_{nd} (1 - \delta_{nd})& n\in \mathcal{N}\setminus{\mathcal{N}_T}, d\in \mathcal{D}_n\\
    \label{eq:eumsp_ref:c3b}   &\phi_{nd} = \sum_{m\in \mathcal{N}_{nd}}\pi_m(r_m^Tx_{m}+\sum_{d'\in \mathcal{D}_m}q_{md'}\delta_{md}+\theta_m)& n\in \mathcal{N}\setminus{\mathcal{N}_{T}}, d\in \mathcal{D}_n\\
    \label{eq:eumsp_ref:c4}    &\theta_{n} = \Theta_n & n\in \mathcal{N}_{T},\\
    \label{eq:eumsp_ref:c5} & x_n \in X_{t(n)} & n\in \mathcal{N},\\
    \label{eq:eumsp_ref:c6} & \delta_{nd} \in \{0,1\} & n\in \mathcal{N},d\in \mathcal{D}_n,\\
    \label{eq:eumsp_ref:c7} & \theta_{n} \in \mathbb{R} & n\in \mathcal{N},\\
    \label{eq:eumsp_ref:c7b} & \phi_{nd} \in \mathbb{R} & n\in \mathcal{N}\setminus{\mathcal{N}_T},d\in \mathcal{D}_n.
  \end{align}
\end{subequations}
Problem \eqref{eq:eumsp_ref} modifies problem \cref{eq:eumsp} in two elements. First, constraints \cref{eq:eumsp_ref:c2} take into
account that $C_{a(n)}=\mathbf{0}$. Second, constraints \cref{eq:eumsp:c3} have been linearized using constants $M_{nd}$ and auxiliary decision variables $\phi_{nd}$
yielding constraints \cref{eq:eumsp_ref:c3,eq:eumsp_ref:c3b} (see e.g., \Cref{sec:app:bigMs} for a general purpose procedure to determine these big-$M$ values).

Problem \cref{eq:eumsp_ref} can be solved in a dynamic programming fashion using the following backward procedure.
The procedure starts by calculating the optimal last-stage expectation for each node at the second-last stage and for
each possible probability distribution.
 For each $m\in \mathcal{N}_{T-1}$ and probability distribution $k\in \mathcal{D}_{m}$,
the optimal expectation at the last stage is obtained by solving the following problem.
\begin{subequations}\label{eq:alg:laststage}
\begin{align}
  \label{eq:alg:laststage_obj}\Phi_{mk} &=\max \sum_{n\in \mathcal{N}_{mk}}\pi_n(r_n^Tx_{n}+\sum_{d\in \mathcal{D}_n}q_{nd}\delta_{nd}+\Theta_n)\\
  \label{eq:alg:laststage_c1}\text{s.t.} & \sum_{d\in\mathcal{D}_n}\delta_{nd} = 1, & n\in \mathcal{N}_{mk},\\
  \label{eq:alg:laststage_c2}&A_nx_n + \sum_{d\in \mathcal{D}_n}B_{nd}\delta_{nd}  = h_n-D_{mk}& n\in \mathcal{N}_{mk},\\
  \label{eq:alg:laststage_c5} & x_n \in X_T & n\in \mathcal{N}_{mk},\\
  \label{eq:alg:laststage_c6} & \delta_{nd} \in \{0,1\} & n\in \mathcal{N}_{mk},d\in \mathcal{D}_n.
\end{align}
\end{subequations}
Problem \eqref{eq:alg:laststage} provides the last-stage expectation for each node $m\in \mathcal{N}_{T-1}$ at the second-last stage, 
and for each possible selection of a probability distribution from $\mathcal{D}_m$. Notice in \cref{eq:alg:laststage} that
$\Theta_n$ is input data and that $D_{mk}$ is moved to the right-hand-side to stress that it does not multiply a decision variable
as in \cref{eq:eumsp_ref:c2},  since the choice of a probability distribution at the parent node $m$ has been fixed to $k$.

 Then, for stage $t = T-2,\ldots,1$, for node $m\in \mathcal{N}_t$, and for distribution $k\in \mathcal{D}_m$ we calculate $\Phi_{km}$ by solving problem \Cref{eq:alg:allstages}.
\begin{subequations}\label{eq:alg:allstages}
\begin{align}
  \label{eq:alg:allstages_obj}\Phi_{mk} &=\max \sum_{n\in \mathcal{N}_{mk}}\pi_n(r_n^Tx_{n}+\sum_{d\in \mathcal{D}_n}q_{nd}\delta_{nd}+\theta_n)\\
  \label{eq:alg:allstages_c1}\text{s.t.} & \sum_{d\in\mathcal{D}_n}\delta_{nd} = 1, & n\in \mathcal{N}_{mk},\\
  \label{eq:alg:allstages_c2}&A_nx_n + \sum_{d\in \mathcal{D}_n}B_{nd}\delta_{nd}  = h_n-D_{mk} & n\in \mathcal{N}_{mk},\\
  \label{eq:alg:allstages_c3}&\theta_n\leq \Phi_{nd} + M_{nd}(1-\delta_{nd}),& d\in \mathcal{D}_n,n\in \mathcal{N}_{mk},\\
  \label{eq:alg:allstages_c4} & x_n \in X_{t(n)} & n\in \mathcal{N}_{mk},\\
  \label{eq:alg:allstages_c5} & \delta_{nd} \in \{0,1\} & n\in \mathcal{N}_{mk},d\in \mathcal{D}_n,\\
  \label{eq:alg:allstages_c6} & \theta_{n} \in \mathbb{R} & n\in \mathcal{N}_{mk}.
\end{align}
\end{subequations}
Notice in problem \cref{eq:alg:allstages} that $\Phi_{nd}$ is input data and has been calculated in the previous steps of the algorithm.

Finally, we can solve the following problem for the root node. 
\begin{subequations}\label{eq:alg:firststage}
\begin{align}
  \label{eq:alg:firststage_obj}z &=\max r_0^Tx_{0}+\sum_{d\in \mathcal{D}_0}q_{0d}\delta_{0d}+\theta_0\\
  \label{eq:alg:firststage_c1}\text{s.t.} & \sum_{d\in\mathcal{D}_0}\delta_{0d} = 1, &\\
  \label{eq:alg:firststage_c2}&A_0x_0 + \sum_{d\in \mathcal{D}_0}B_{0d}\delta_{0d}  = h_0, & \\
  \label{eq:alg:firststage_c3}&\theta_0\leq \Phi_{0d} + M_{0d}(1-\delta_{0d}),& d\in \mathcal{D}_0,\\
  \label{eq:alg:firststage_c4} & x_0 \in X_{1}, & \\
  \label{eq:alg:firststage_c5} & \delta_{0d} \in \{0,1\} & d\in \mathcal{D}_0,\\
  \label{eq:alg:firststage_c6} & \theta_{0} \in \mathbb{R}. &
\end{align}
\end{subequations}

Letting $D=\max_{n\in\mc{N}}|\mc{D}_n|$, the algorithm entails solving $\mc{O}\big(D|\mc{N}|\big)$ mixed-integer programs of size significantly smaller than \eqref{eq:eumsp}.

\section{Case study}\label{sec:cs}

In this section we present a case study based on the \textit{Football Team Composition Problem} (FTCP -- \citet{Pan17,PanH20})
which we extend in order to account for decision-dependent uncertainty. The problem consists
of selecting players for a football team while their future market value is uncertain and influenced
by the team for which they play. The scope of the computational study is to compare the performance of the algorithm
to that of a state-of-the-art commercial solver. 
We remark, however, that formulation \cref{eq:eumsp} and the solution algorithm proposed are general and applicable beyond the context of the FTCP,
which we use solely as an example. Decision problems under endogenous uncertainty may, in fact, arise in several business context, some of which are mentioned in \Cref{sec:introduction}.
To provide another practical example, consider an agribusiness making periodic production planning decisions for a number of different crops while demand is uncertain.
Product substitution is a common phenomenon in the production of crops. In fact, a customer, say a farmer, may view multiple crops as suitable for their farming needs,
see, e.g., \citet{BanD20}. Thus, the uncertain demand for a crop may depend in part on the portfolio of crops an agribusiness offers for sale, yielding a decision problem under endogenous uncertainty.

In \Cref{sec:cs:ftcp} we describe the FTCP in more details and formulate it as
an MSPEU. In \Cref{sec:cs:bigM} we illustrate an efficient procedure to obtain big-$M$ values.
In \Cref{sec:cs:instances} we introduce the problem instances and finally in \Cref{sec:cs:results} we
present and discuss the results. The data of the problem instances is available online at \url{https://github.com/GioPan/instancesFTCPwithEndogenousUncertainty}
for the benefit of future research. 

\subsection{The Football Team Composition Problem}\label{sec:cs:ftcp}
The FTCP is the problem of composing a football team by purchasing and selling professional football players.
A complete description of the problem can be found in \citet{Pan17}. In what follows, we
report the basic elements necessary for this case study. The decision problem can be described as follows.
At every stage, i.e., \textit{transfer market window} (TMW), professional football clubs can renew their teams
by purchasing players from other clubs or selling available players to other clubs.
In order to participate in national and international competitions, professional clubs must compose a team made
of a fixed number of players. In addition, the coach of the team requires players with different roles (e.g., defenders or mid-fielders) and
skills. 

Club managers are often given a budget to spend in the transfer market and are typically allowed to
reinvest the revenue from the sale of players. The current market value of football players is known, while the
future value is stochastic as it depends on a number of random events such as injuries,
fitness, motivation and ultimately luck. Therefore, at every TMW, football clubs make decisions in conditions of uncertainty
with the scope of maximizing the expected value of the team. Furthermore, the market values of the players in the same
team are strongly correlated. This generates a multistage stochastic program with endogenous uncertainty, since the decision of
hiring or selling a player will change the correlations of the joint value distribution of the players considered.
That is, hiring a football player will make their value strongly correlated with the value of the other players in the team, while selling or not hiring a
player will make their value uncorrelated with the value of the players in the team. 

In order to model the FTCP we assume the club is evaluating a finite number of alternative \textit{team compositions}. Every team composition
consists of the required number of players and contains the necessary mix of skills. Let $\mathcal{I}$ be the set of possible team compositions,
$\mathcal{N}$ the set of nodes in the MDST describing the underlying uncertainty, $\mathcal{N}_{in}$
the set of child nodes of node $n\in \mathcal{N}$ if team composition $i\in\mc{I}$ is chosen at node $n$.
Note in fact that the selection of a team composition determines the correlations between the values of the players in the instance, and thus
the probability distribution. Let $\mathcal{N}^L$ be the set of leaf nodes,
i.e., the nodes at the last decision stage, and $t(n)\in \{1,\ldots,T\}$ the stage corresponding to node $n$. 

Let $V_{in}$ be the value of team composition $i\in \mc{I}$ at node $n$ (i.e., the sum of the values of the players in team
composition $i$), let $C^S_{in}$ be the cost of the salary of the players in team composition $i$ at node $n$, $C^T_{jn}$
the cost of transitioning from team composition $i$ to team composition $j$ at node $n$, i.e., the cost of the transfer fees for the players bought
minus the revenue for the transfer fees of the players sold.  Furthermore, let $B$ be the budget available to the club for the transfer market at node $n$
and $E_n$ a stochastic extra-budget conditional on events such as sport successes. 
Let $X_i$ be equal to one if team composition $i\in \mathcal{I}$ is the initial team composition, $M_{in}$ a suitably large constant for each $i$ and $n$ (see \Cref{sec:cs:bigM}),
and $\rho$ a discount rate. Let $\delta_{in}$ be a binary decision variable which is equal $1$ if team composition $i\in \mathcal{I}$ is chosen at node
$n\in \mathcal{N}$, $x_{ijn}$ a binary decision variable which is equal $1$ if the club transitions from team composition $i\in \mathcal{I}$ to team
composition $j\in \mathcal{J}$ at node $n\in \mathcal{N}$ and, $\phi_{in}$ the expected net value of the team at the children of node $n$ if composition $i$ is chosen and, finally,
$\theta_n$ the expected net value of the team at the child nodes of node $n\in \mathcal{N}$.
The FTCP with endogenous uncertainty is hence:

\begin{subequations}\label{ftcp2}
  \begin{align}
   \label{ftcp2:obj} \max & \sum_{i\in \mathcal{I}}(V_{i0}\delta_{i0}-C^S_{i0}\delta_{i0} - \sum_{j\in \mathcal{I}}C^T_{ij0}x_{ij0})+\theta_0 & \\
   \label{ftcp2:onlyone} \text{s.t.}& \sum_{i\in\mathcal{I}}\delta_{in} = 1 & n\in \mathcal{N},\\
   \label{ftcp2:changeorhold0} &\delta_{i0} - X_{i} + \sum_{j\in \mathcal{I}:j\neq i}(x_{ij0} - x_{ji0}) = 0 &i\in \mathcal{I},\\
    \label{ftcp2:changeorhold}&\delta_{in} - \delta_{i,a(n)} + \sum_{j\in \mathcal{I}:j\neq i}(x_{ijn} - x_{jin}) = 0 &n\in \mathcal{N}\setminus{\{0\}},i\in \mathcal{I},\\
    \label{ftcp2:maxOneChange}&\sum_{i\in\mc{I}}\sum_{j\neq i \in\mc{I}}x_{ijn} \leq 1 &n\in \mathcal{N},\\
    \label{ftcp2:budget}&\sum_{i\in \mathcal{I}}\sum_{j\in \mathcal{I}}C^T_{ijn}x_{ijn}\leq B+E_n &n\in \mathcal{N},\\
    \label{ftcp2:theta}&\theta_n\leq \phi_{in}+M_{in}(1-\delta_{in})& i\in \mathcal{I},n\in \mathcal{N},\\
    \label{ftcp2:phi}&\phi_{in} =  \sum_{m\in \mathcal{N}_{in}}\frac{1}{1+\rho^{t(m)}}\pi_m\left\{\theta_m + \sum_{j\in \mathcal{I}}\left(\vphantom{\sum_{i\in C}}V_{jm}\delta_{jm}-C^S_{jm}\delta_{jm}-\sum_{k\in \mathcal{I}}C^T_{jkm}x_{jkm}\right)\right\}& i\in \mathcal{I},n\in \mathcal{N}\setminus{\mathcal{N}^L},\\
   \label{ftcp2:thetafinal}&\theta_{n} = \frac{1}{1+\rho^{t(n)+1}}\sum_{i\in \mathcal{I}}V_{in}\delta_{in} & n\in \mathcal{N}^L,\\
    \label{ftcp2:range1}&\delta_{in}\in \{0,1\}& i\in \mathcal{I},n\in \mathcal{N},\\
    \label{ftcp2:range2}&x_{ijn}\in \{0,1\}& i,j\in \mathcal{I},n\in \mathcal{N},\\
    \label{ftcp2:range3}&\theta_{n} \in \mathbb{R}& n\in \mathcal{N}.
  \end{align}
\end{subequations}

Objective function \eqref{ftcp2:obj} represents the sum of the value of the team composition chosen here-and-now, minus the expenses for salaries and, if any, the transition cost from the initial team composition. In addition it takes into account the expected value of the team at future nodes. Constraints \eqref{ftcp2:onlyone} ensure that only one team composition is chosen at each decision node, while constraints \eqref{ftcp2:changeorhold0}-\eqref{ftcp2:changeorhold} state that either the club holds the same team composition as in the previous season or a new team composition is chosen at node $0$ and at the rest of the nodes, respectively.
Constraints \eqref{ftcp2:maxOneChange} state that at most one team change can be done at each node.
Constraints \eqref{ftcp2:budget} ensure that the net expenses for transitioning from a team composition to another (i.e., the money spent in the transfer market)
do not exceed the available budget.
Constraints \eqref{ftcp2:theta} state that the future net expected value of the team depends on the team composition chosen. Constraints \eqref{ftcp2:phi} sets the expected value at the children of node $n$ if composition $i$ is chosen, while constraints \eqref{ftcp2:thetafinal} set the final value of $\theta_n$ to the value of the team chosen at the leaf nodes. This corresponds to a sunset value which accounts the termination of an infinite horizon problem. Notice that the final value of the team is discounted as a future value. Finally, constraints \eqref{ftcp2:range1}-\eqref{ftcp2:range3} define the domain for the decision variables.

\subsection{Finding big-M values for the FTCP}\label{sec:cs:bigM} 
We illustrate a fast method to set the $M_{in}$ values in \eqref{ftcp2}. The method proposed resulted, on this special case, significantly faster
than the general method described in \Cref{sec:app:bigMs}.

We start by observing that, for all $i\in \mathcal{I}$ and $n\in \mathcal{N}$ we are looking for a value $M_{in}$ such that:
$$\theta_n - \phi_{in}\leq M_{in}$$
and that, based on constraints \cref{ftcp2:theta}, we have
$$\theta_n \leq \max_{i\in \mc{I}} \phi_{in}$$
consequently, we need to set 
$$M_{in}\geq \max_{i\in \mc{I}} \phi_{in} - \phi_{in}$$
thus, under the very mild assumption that $\phi_{in}\geq 0$ (i.e., that the team does not spend in transfers and salaries more than the value of the entire team), for every $i\in \mc{I}$ a valid $M_{in}$ is
$$M_{in} = M_n = \max_{i\in \mc{I}} \phi_{in}$$
By expanding $\phi_{in}$ according to constraints \cref{ftcp2:phi} and \cref{ftcp2:thetafinal} we can set
\begin{itemize}
\item $M_{in} = \max_{i\in \mc{I}}\frac{1}{1+\rho^{t(n)+1}}V_{in}$ for all $i\in\mc{I}$ if $n\in \mathcal{N}^L$
\item $M_{in} = \max_{i\in \mc{I}} \Bigg\{ \sum_{m\in \mathcal{N}_{in}}\frac{1}{1+\rho^{t(m)}}\pi_m\bigg( M_m + \max_{j\in \mathcal{I}}\left\{\vphantom{\sum_{i\in C}} V_{jm}\right\}\bigg)\Bigg\}$ for all $i\in \mc{I}$ if $n\in \mc{N}\setminus{\mc{N}^L}$.
\end{itemize}

\subsection{Problem instances}\label{sec:cs:instances}

Instances are generated from the case studies on the FTCP presented by \citet{Pan17} and are available online at \url{https://github.com/GioPan/instancesFTCPwithEndogenousUncertainty}.
The instances consist of the $20$ teams competing in the English Premier League 2013/2014 and dealing
with the transfer market of summer 2014. Each team is characterized by a budget, a list of players (made of the current and target players), and a number of randomly generated team compositions
$\mc{I}$ complying with the regulations of the league and with the coach's specifics. In turn, each player is characterized by age, role, current market value, salary, selling and purchase price.
The value and cost of salaries of a team composition are calculated as the sum of values and salaries, respectively, of the players contained. Similarly, transition costs are calculated as the revenue for the players sold, minus the cost of the players bought.

Future player values are stochastic and modeled by means of the regression equation available in
\citet{Pan17}. The joint probability distribution of the market values for all the players forms a
multivariate normal distribution. In addition, we set a $0.8$ correlation between the players belonging to the focal team
and no correlation between the players belonging to the focal team and the remaining players. Therefore, the correlations
change with the decisions of the team to buy or sell players creating an MSPEU. 
Scenario trees are obtained by sampling realizations from the underlying multivariate normal distributions.
 Finally, for each node in the scenario tree, the random extra-budget $E_n$ is calculated as a percentage of the deterministic budget $B$.
The percentage is the same as percentage increase of team value from the parent node, that is $\max\left\{0,100 * (V_{in}-V_{i,a(n)})/V_{i,a(n)}\right\}$.
That is, an increase in the team value yields an increase in the spending capability of the team. 

The instances represent a four-stage horizon and are identified by (i) the focal team, where Team $\in\{$ ARS, ASV, CAR, CHE, CRP, EVE, FUL, HUL, LIV, MAC, MAU, NEC, NOR, SOU, STO, SUN, SWA, TOT, WBA, WHU $\}$, (ii) the number of team composition $|\mc{I}|\in\{3,4,5\}$, and (iii) the number of samples in each distribution $S\in \{4,5\}$, yielding 120 numerically different instances with up to approximately $400$ thousand binary variables.

\subsection{Results}\label{sec:cs:results}
We implemented our algorithm in Python 2.7 using Cplex 12.8 for solving the subproblems. Cplex 12.8 has also been used to solve the full problems. 
All experiments have been run on a server with 64 double AMD Opteron 6380 processors and 251 GB RAM. 
\Cref{tab:resC3,tab:resC4,tab:resC5} report the results for different cardinalities $|\mc{I}|$ and different number of samples $S$ for all the focal teams.
The computation times do not include the time required to find big-$M$ values as illustrated in \Cref{sec:cs:bigM}.
An account of these times is reported in \Cref{app:bigMtime}. Note that big-$M$ values are required both when solving the full problem and when using our algorithm and in both cases are calculated beforehand. Thus, their computation time does not have an impact on the comparison between the two solution methods. Since the algorithm entails iteratively building and solving mixed-integer programs, the times reported include the time required for building the models. 

\Cref{tab:resC3} reports the results with $|\mc{I}|=3$, generating instances with up to approximately 32000 binary variables.
For these instances it can be noticed that Cplex performs much better than our method on the smaller test cases (with $S=4$ samples)
while our algorithm is more competitive on the instances with $S=5$. Altogether, Cplex solves the problems approximately 25\% faster than our method.

\begin{longtable}{lc|ccc|cccc|c}
  \caption{Results with $|\mc{I}| = 3$. $S$ indicates the number of realizations describing each distribution. \#Var, \#Bin and \#Con indicate the total number of variables, the number of binary variables, and the number of constraints, respectively. $\tau$ (Cpx) indicate the elapsed time when using the algorithm (Cplex). Obj. (Cpx) indicates the objective value obtained using the algorithm (Cplex). $\Delta \tau$ is calculated as $100$($\tau$ - $\tau$ Cpx)/$\tau$ Cpx.}\label{tab:resC3}\\
  \toprule
Team &    $S$ &  \#Var &    \#Bin &  \#Con &   $\tau$ Cpx [sec]&  Obj. Cpx &  $\tau$ [sec]&  Obj.&  $\Delta \tau$ [\%]\\
\midrule
 ARS &  4 &       24505 &      16965 &         22620 &    86.924 &     1487.985 &  178.413 &     1487.985 &    105.252 \\
 ASV &  4 &       24505 &      16965 &         22620 &   107.462 &       858.189 &  228.801 &       858.189 &    112.914 \\
 CAR &  4 &       24505 &      16965 &         22620 &   112.806 &       642.253 &  282.521 &       642.253 &    150.448 \\
 CHE &  4 &       24505 &      16965 &         22620 &   144.798 &     2604.929 &  297.918 &     2604.929 &    105.748 \\
 CRP &  4 &       24505 &      16965 &         22620 &   105.784 &       471.838 &  232.649 &       471.838 &    119.929 \\
 EVE &  4 &       24505 &      16965 &         22620 &   198.908 &     1216.299 &  189.362 &     1216.299 &     -4.799 \\
 FUL &  4 &       24505 &      16965 &         22620 &   116.068 &       694.942 &  187.200 &       694.942 &     61.284 \\
 HUL &  4 &       24505 &      16965 &         22620 &   174.406 &       571.997 &  215.792 &       571.997 &     23.729 \\
 LIV &  4 &       24505 &      16965 &         22620 &   457.376 &     1371.466 &  196.851 &     1371.466 &    -56.961 \\
 MAC &  4 &       24505 &      16965 &         22620 &   165.375 &     1734.123 &  166.470 &     1734.123 &      0.662 \\
 MAU &  4 &       24505 &      16965 &         22620 & 2036.260 &     2474.198 &  231.169 &     2474.198 &    -88.647 \\
 NEC &  4 &       24505 &      16965 &         22620 &    95.849 &       974.510 &  213.350 &       974.510 &    122.589 \\
 NOR &  4 &       24505 &      16965 &         22620 &    69.662 &       486.204 &  161.434 &       486.204 &    131.739 \\
 SOU &  4 &       24505 &      16965 &         22620 &   153.804 &       640.762 &  143.529 &       640.762 &     -6.681 \\
 STO &  4 &       24505 &      16965 &         22620 &    79.237 &       674.964 &  165.104 &       674.964 &    108.366 \\
 SUN &  4 &       24505 &      16965 &         22620 &   229.976 &       837.582 &  202.287 &       837.582 &    -12.040 \\
 SWA &  4 &       24505 &      16965 &         22620 &    90.615 &       659.840 &  228.751 &       659.840 &    152.442 \\
 TOT &  4 &       24505 &      16965 &         22620 &   193.830 &     1383.202 &  171.743 &     1383.202 &    -11.395 \\
 WBA &  4 &       24505 &      16965 &         22620 &   223.691 &       570.527 &  179.085 &       570.527 &    -19.941 \\
 WHU &  4 &       24505 &      16965 &         22620 &   259.587 &       785.027 &  202.723 &       785.027 &    -21.905 \\
 ARS &  5 &       47008 &      32544 &         43392 &   189.159 &     1508.956 &  403.177 &     1508.956 &    113.142 \\
 ASV &  5 &       47008 &      32544 &         43392 &   486.235 &       801.662 &  400.098 &       801.662 &    -17.715 \\
 CAR &  5 &       47008 &      32544 &         43392 &   223.908 &       661.053 &  389.787 &       661.053 &     74.083 \\
 CHE &  5 &       47008 &      32544 &         43392 &   593.065 &     2667.682 &  628.685 &     2667.682 &      6.006 \\
 CRP &  5 &       47008 &      32544 &         43392 &   222.443 &       478.977 &  483.007 &       478.977 &    117.138 \\
 EVE &  5 &       47008 &      32544 &         43392 &   726.302 &     1134.583 &  387.368 &     1134.583 &    -46.666 \\
 FUL &  5 &       47008 &      32544 &         43392 &   275.205 &       717.025 &  377.999 &       717.025 &     37.352 \\
 HUL &  5 &       47008 &      32544 &         43392 & 2785.695 &       494.981 &  446.281 &       494.981 &    -83.980 \\
 LIV &  5 &       47008 &      32544 &         43392 & 7386.953 &     1403.374 &  391.775 &     1403.374 &    -94.696 \\
 MAC &  5 &       47008 &      32544 &         43392 &   349.382 &     1769.971 &  306.545 &     1769.971 &    -12.261 \\
 MAU &  5 &       47008 &      32544 &         43392 & 7755.757 &     2260.280 &  488.009 &     2260.280 &    -93.708 \\
 NEC &  5 &       47008 &      32544 &         43392 &   481.964 &     1014.031 &  440.037 &     1014.031 &     -8.699 \\
 NOR &  5 &       47008 &      32544 &         43392 &   356.411 &       499.601 &  317.165 &       499.601 &    -11.012 \\
 SOU &  5 &       47008 &      32544 &         43392 &   148.499 &       656.613 &  322.197 &       656.613 &    116.968 \\
 STO &  5 &       47008 &      32544 &         43392 &   391.803 &       689.878 &  359.498 &       689.878 &     -8.245 \\
 SUN &  5 &       47008 &      32544 &         43392 &   481.092 &       854.091 &  406.024 &       854.091 &    -15.604 \\
 SWA &  5 &       47008 &      32544 &         43392 &   444.930 &       681.336 &  386.791&       681.336 &    -13.067 \\
 TOT &  5 &       47008 &      32544 &         43392 &   418.714 &     1409.204 &  400.087&     1409.204 &     -4.449 \\
 WBA &  5 &       47008 &      32544 &         43392 &   299.717 &       522.849 &  367.367 &       522.849 &     22.571 \\
 WHU &  5 &       47008 &      32544 &         43392 &   654.066 &       678.902 &  426.684 &       678.902 &    -34.764 \\
 \midrule
Avg&             & &   &  & 744.342&        &305.093  &   &      25.378 \\
 \bottomrule
\end{longtable}

By increasing the number of team compositions we generate larger MDSTs, and the corresponding MSPEUs also become larger. \Cref{tab:resC4} reports the results with $|\mc{I}|=4$ generating a MDST with four possible distributions at each decision node and problems with up to approximately 135 thousand binary variables. In this case, our algorithm performs significantly better that Cplex on almost all instances, solving the problems 32.8\% faster.

\begin{longtable}{lc|ccc|cccc|c}
  \caption{Results with $|\mc{I}| = 4$. $S$ indicates the number of realizations describing each distribution. \#Var, \#Bin and \#Con indicate the total number of variables, the number of binary variables, and the number of constraints, respectively. $\tau$ (Cpx) indicate the elapsed time when using the algorithm (Cplex). Obj. (Cpx) indicates the objective value obtained using the algorithm (Cplex). $\Delta \tau$ is calculated as $100$($\tau$ - $\tau$ Cpx)/$\tau$ Cpx.}\label{tab:resC4}\\
  \toprule
Case &   $S$ &  \#Var &    \#Bin &  \#Con &   $\tau$ Cpx [sec]&  Obj. Cpx &  $\tau$ [sec]&  Obj.&  $\Delta \tau$ [\%]\\
\midrule
 ARS &4 &       91749 &    69904 &         65535 &   469.440 &     1409.074 &   585.920 &     1409.074 &     24.812 \\
 ASV &4 &       91749 &    69904 &         65535 &   560.788 &       916.248 &   453.329 &       916.248 &    -19.162 \\
 CAR &4 &       91749 &    69904 &         65535 &   551.781 &       574.859 &   524.636 &       574.859 &     -4.920 \\
 CHE &4 &       91749 &    69904 &         65535 &   756.516 &     2476.158 &   729.832 &     2476.158 &     -3.527 \\
 CRP &4 &       91749 &    69904 &         65535 &   558.236 &       431.065 &   632.984 &       431.065 &     13.390 \\
 EVE &4 &       91749 &    69904 &         65535 & 7656.373 &     1282.401 &   535.382 &     1282.401 &    -93.007 \\
 FUL &4 &       91749 &    69904 &         65535 &   932.283 &       711.560 &   488.729 &       711.560 &    -47.577 \\
 HUL &4 &       91749 &    69904 &         65535 & 7746.679 &       559.913 &   608.927 &       559.913 &    -92.140 \\
 LIV &4 &       91749 &    69904 &         65535 & 7627.856 &     1452.961 &   545.440 &     1452.961 &    -92.849 \\
 MAC &4 &       91749 &    69904 &         65535 &   888.410 &     1657.326 &   371.674 &     1657.326 &    -58.164 \\
 MAU &4 &       91749 &    69904 &         65535 & 7672.755 &     2489.203 &   701.230 &     2489.203 &    -90.861 \\
 NEC &4 &       91749 &    69904 &         65535 &   545.592 &       908.271 &   518.564 &       908.271 &     -4.954 \\
 NOR &4 &       91749 &    69904 &         65535 &   393.393 &       443.482 &   371.138 &       443.482 &     -5.657 \\
 SOU &4 &       91749 &    69904 &         65535 &   368.125 &       603.976 &   342.766 &       603.976 &     -6.889 \\
 STO &4 &       91749 &    69904 &         65535 &   453.101 &       633.148 &   396.288 &       633.148 &    -12.539 \\
 SUN &4 &       91749 &    69904 &         65535 &   529.445 &       785.688 &   513.954 &       785.688 &     -2.926 \\
 SWA &4 &       91749 &    69904 &         65535 &   469.695 &       607.176 &   440.718 &       607.176 &     -6.169 \\
 TOT &4 &       91749 &    69904 &         65535 &   639.796 &     1295.453 &   421.155 &     1295.453 &    -34.174 \\
 WBA &4 &       91749 &    69904 &         65535 & 3945.793 &       633.948 &   473.058 &       633.948 &    -88.011 \\
 WHU &4 &       91749 &    69904 &         65535 & 8498.951 &       828.664 &   601.910 &       828.664 &    -92.918 \\
 ARS &5 &      176841 &   134736 &        126315 & 1210.558 &     1440.502 & 1340.235 &     1440.502 &     10.712 \\
 ASV &5 &      176841 &   134736 &        126315 & 1946.890 &       760.043 & 1381.227 &       760.043 &    -29.055 \\
 CAR &5 &      176841 &   134736 &        126315 & 1357.787 &       597.870 & 1329.530 &       597.870 &     -2.081 \\
 CHE &5 &      176841 &   134736 &        126315 & 1706.027 &     2602.808 & 1896.954 &     2602.808 &     11.191 \\
 CRP &5 &      176841 &   134736 &        126315 & 1365.008 &       467.102 & 1489.021 &       467.102 &      9.085 \\
 EVE &5 &      176841 &   134736 &        126315 & 8505.074 &     1232.924 & 1260.209 &     1232.924 &    -85.183 \\
 FUL &5 &      176841 &   134736 &        126315 & 7808.254 &       699.597 & 1388.999 &       699.597 &    -82.211 \\
 HUL &5 &      176841 &   134736 &        126315 & 8642.826 &       553.807 & 1300.262 &       553.807 &    -84.956 \\
 LIV &5 &      176841 &   134736 &        126315 & 8326.633 &     1460.389 & 1234.644 &     1460.389 &    -85.172 \\
 MAC &5 &      176841 &   134736 &        126315 & 1020.989 &     1678.841 & 1149.133 &     1678.841 &     12.551 \\
 MAU &5 &      176841 &   134736 &        126315 & 8397.024 &     2285.317 & 1405.108 &     2285.317 &    -83.267 \\
 NEC &5 &      176841 &   134736 &        126315 & 1623.349 &     1062.903 & 1427.509 &     1062.903 &    -12.064 \\
 NOR &5 &      176841 &   134736 &        126315 & 1041.199 &       480.087 & 1085.893 &       480.087 &      4.293 \\
 SOU &5 &      176841 &   134736 &        126315 & 1000.712 &       616.957 & 1020.687 &       616.957 &      1.996 \\
 STO &5 &      176841 &   134736 &        126315 & 1111.253 &       663.486 & 1059.614 &       663.486 &     -4.647 \\
 SUN &5 &      176841 &   134736 &        126315 & 1480.353 &       891.809 & 1609.855 &       891.809 &      8.748 \\
 SWA &5 &      176841 &   134736 &        126315 & 1197.202 &       623.004 & 1596.017 &       623.004 &     33.312 \\
 TOT &5 &      176841 &   134736 &        126315 & 2560.338 &     1441.642 & 1264.625 &     1441.642 &    -50.607 \\
 WBA &5 &      176841 &   134736 &        126315 & 8304.973 &       601.984 & 1152.420 &       601.984 &    -86.124 \\
 WHU &5 &      176841 &   134736 &        126315 & 8432.614 &       755.708 & 1416.893 &       755.708 &    -83.197 \\
\midrule
Avg&             & &   &  & 3207.601&        & 926.661 &   &      -32.873 \\
 \bottomrule
\end{longtable}

This pattern illustrates that our algorithm scales significantly better than Cplex, and is confirmed also by our results on the instances with $|\mc{I}|=5$, the largest we tested. These instances generate problems with up to approximately 400 thousand binary variables and the corresponding results are provided in \Cref{tab:resC5}.
 It is possible to notice that our algorithm outperforms Cplex on almost all instances. Altogether, our algorithm is able to solve the same instances 61.1\% faster than Cplex.
Only on two instances (CAR and CHE, with $S=4$) Cplex performs better than our algorithm.
Similar cases are to be expected since our algorithm entails solving several subproblems which are mixed-integer programs.
Therefore, it is possible that some numerically difficult subproblems are found that slow down the entire algorithm (in our case we do not solve problems in parallel).
However, despite that, the average performance of our algorithm on the largest instance is by far better.

These results illustrate that MSPEU easily generate very large optimization problems. The size of the problems, in our instances, increases by approximately three to five times by increasing the number of possible distributions.
However, the algorithm we provide scales better than the solver Cplex. Particularly, when using Cplex, the average solution time increases by approximately 330\% when increasing $|\mc{I}|$ from 3 to 4, and by approximately 175\% when increasing $|\mc{I}|$ from 4 to 5. With our algorithm the average solution time increases by approximately 203\% when increasing $|\mc{I}|$ from 3 to 4, and by approximately 158\% when increasing $|\mc{I}|$ from 4 to 5.

\begin{longtable}{lc|ccc|cccc|c}
  \caption{Results with $|\mc{I}| = 5$. $S$ indicates the number of realizations describing each distribution. \#Var, \#Bin and \#Con indicate the total number of variables, the number of binary variables, and the number of constraints, respectively. $\tau$ (Cpx) indicate the elapsed time when using the algorithm (Cplex). Obj. (Cpx) indicates the objective value obtained using the algorithm (Cplex). $\Delta \tau$ is calculated as $100$($\tau$ - $\tau$ Cpx)/$\tau$ Cpx.}\label{tab:resC5}\\
  \toprule
Case & $S$ &  \#Var &    \#Bin &  \#Con &   $\tau$ Cpx [sec]&  Obj. Cpx &  $\tau$ [sec]&  Obj.&  $\Delta \tau$ [\%]\\
\midrule
 ARS &4 &      261051 &      210525 &        151578 &  4148.394 &     1412.322 & 2487.630 &     1412.322 &    -40.034 \\
 ASV &4 &      261051 &      210525 &        151578 &  9128.747 &       888.394 & 2360.369 &       888.394 &    -74.144 \\
 CAR &4 &      261051 &      210525 &        151578 &  2324.632 &       559.909 & 2951.872 &       559.909 &     26.982 \\
 CHE &4 &      261051 &      210525 &        151578 &  2941.786 &     2351.709 & 4074.808 &     2351.709 &     38.515 \\
 CRP &4 &      261051 &      210525 &        151578 &  4873.336 &       403.008 & 2757.925 &       403.008 &    -43.408 \\
 EVE &4 &      261051 &      210525 &        151578 &  9240.398 &     1318.359 & 2367.225 &     1328.380 &    -74.382 \\
 FUL &4 &      261051 &      210525 &        151578 &  9385.279 &       676.308 & 2775.651 &       676.308 &    -70.425 \\
 HUL &4 &      261051 &      210525 &        151578 &  9616.179 &       496.264 & 2430.993 &       600.252 &    -74.720 \\
 LIV &4 &      261051 &      210525 &        151578 &  9892.436 &     1584.601 & 1522.159 &     1584.601 &    -84.613 \\
 MAC &4 &      261051 &      210525 &        151578 &  1748.831 &     1589.133 & 1070.475 &     1589.133 &    -38.789 \\
 MAU &4 &      261051 &      210525 &        151578 &  9294.710 &     2529.355 & 1935.253 &     2533.948 &    -79.179 \\
 NEC &4 &      261051 &      210525 &        151578 &  6984.354 &     1202.156 & 1502.823 &     1202.156 &    -78.483 \\
 NOR &4 &      261051 &      210525 &        151578 &  3516.805 &       420.312 & 1038.191 &       420.312 &    -70.479 \\
 SOU &4 &      261051 &      210525 &        151578 &  3593.291 &       599.703 &   947.766 &       599.703 &    -73.624 \\
 STO &4 &      261051 &      210525 &        151578 &  1966.818 &       632.875 & 1182.410 &       632.875 &    -39.882 \\
 SUN &4 &      261051 &      210525 &        151578 &  4767.951 &       731.767 & 1460.145 &       731.767 &    -69.376 \\
 SWA &4 &      261051 &      210525 &        151578 &  4397.444 &       576.314 & 1280.492 &       576.314 &    -70.881 \\
 TOT &4 &      261051 &      210525 &        151578 &  3894.520 &     1392.754 & 1169.367 &     1392.754 &    -69.974 \\
 WBA &4 &      261051 &      210525 &        151578 & 11192.050 &       590.701 & 1298.701 &       661.566 &    -88.396 \\
 WHU &4 &      261051 &      210525 &        151578 &  9314.423 &       600.960 & 1755.932 &       768.724 &    -81.148 \\
 ARS &5 &      504556 &      406900 &        292968 &  5985.321 &     1453.229 & 3231.294 &     1453.229 &    -46.013 \\
 ASV &5 &      504556 &      406900 &        292968 & 13100.017 &       969.935 & 2557.610 &       969.935 &    -80.476 \\
 CAR &5 &      504556 &      406900 &        292968 & 12990.803 &       612.171 & 3115.356 &       612.171 &    -76.019 \\
 CHE &5 &      504556 &      406900 &        292968 & 15815.020 &     2608.034 & 4221.096 &     2608.034 &    -73.310 \\
 CRP &5 &      504556 &      406900 &        292968 & 13160.004 &       476.955 & 3668.813 &       476.955 &    -72.121 \\
 EVE &5 &      504556 &      406900 &        292968 & 18368.885 &     1289.379 & 2913.376 &     1452.827 &    -84.140 \\
 FUL &5 &      504556 &      406900 &        292968 & 13385.821 &       746.554 & 2811.966 &       755.759 &    -78.993 \\
 HUL &5 &      504556 &      406900 &        292968 & 12995.828 &       516.157 & 3335.790 &       655.447 &    -74.332 \\
 LIV &5 &      504556 &      406900 &        292968 & 12808.767 &     1473.064 & 2854.447 &     1706.759 &    -77.715 \\
 MAC &5 &      504556 &      406900 &        292968 &  5341.258 &     1681.019 & 2050.793 &     1681.019 &    -61.605 \\
 MAU &5 &      504556 &      406900 &        292968 & 13282.279 &     2690.920 & 3722.067 &     2712.864 &    -71.977 \\
 NEC &5 &      504556 &      406900 &        292968 & 13368.708 &     1365.929 & 2971.268 &     1365.929 &    -77.774 \\
 NOR &5 &      504556 &      406900 &        292968 &  5352.415 &       462.292 & 2111.632 &       462.292 &    -60.548 \\
 SOU &5 &      504556 &      406900 &        292968 &  5347.118 &       624.205 & 1897.896 &       624.205 &    -64.506 \\
 STO &5 &      504556 &      406900 &        292968 &  5679.705 &       660.999 & 2308.348 &       660.999 &    -59.358 \\
 SUN &5 &      504556 &      406900 &        292968 & 13483.620 &       927.387 & 2961.979 &       927.387 &    -78.033 \\
 SWA &5 &      504556 &      406900 &        292968 & 11568.065 &       631.495 & 2500.653 &       631.495 &    -78.383 \\
 TOT &5 &      504556 &      406900 &        292968 & 12955.368 &     1490.643 & 2323.134 &     1490.643 &    -82.068 \\
 WBA &5 &      504556 &      406900 &        292968 & 12852.950 &       709.452 & 2696.335 &       718.245 &    -79.022 \\
 WHU &5 &      504556 &      406900 &        292968 & 13408.572 &       843.851 & 3371.308 &       844.027 &    -74.857 \\
\midrule
Avg&             & &   &  & 8836.822&        &  2399.883&   &      -65.192 \\
 \bottomrule
\end{longtable}

\section{Conclusions}\label{sec:conclusions}
This paper introduced 1) a novel scenario tree structure and 2) a node formulation for multistage stochastic programs with endogenous uncertainty, as well as 3) a solution algorithm for a special case.
A computational study shows that while, as expected, problems with endogenous uncertainty tend to generate large optimization problems, all our instances where solvable by Cplex to optimality in at most approximately 5 hours. Furthermore, our algorithm outperformed Cplex on the medium and large instances and showed that it scales well with the size of the problem.

Despite the encouraging results obtained in our study, solving multistage stochastic programs with endogenous uncertainty remains, in general, a challenging task.
Our algorithm requires an explicit scenario tree structure, and solves a number of problems which grows linearly with the number of nodes.
However, the number of nodes in a scenario tree grows exponentially with the number of stages and the treatment of cases with more than a handful of stages may soon become prohibitive.
New approaches in the spirit of \citet{PerP91,ZouAS19}, based on progressive approximations of future stages, may be proven more scalable.
Furthermore, our models employ so called big-$M$ constants. Poorly chosen big-$M$ values, e.g., by trial-and-error, may become problematic. It is well known that they may create numerical difficulties when solving mixed-integer programs. In addition, as discussed in \citet{PinM19}, they may lead to highly sub-optimal solutions. The authors use a simple bilevel programming problem (which can be reformulated as a mixed-integer program that uses big-$M$s) to show how a poorly designed trial-and-error procedure may generate the false belief that the solution to the reformulation is indeed optimal for the original bilevel problem. Based on these evidences we also advocate caution and the use of more sophisticated procedures for setting big-$M$ values. The procedure in \Cref{sec:app:bigMs} goes in this direction.
Finally, cases more general than the special one treated \Cref{sec:algorithm} remain to be addressed.

\bibliography{ftcp_endogenous_uncertainty}

\begin{appendices}
  \crefalias{section}{appsec}
  \section{Notation table}\label{app:notation}
  
  \begin{longtable}{p{0.22\textwidth} p{0.78\textwidth}}
    \caption{Notation of problem \eqref{eq:eumsp}.}\\
    \toprule
    \multicolumn{2}{c}{Sets}  \\
    \midrule
    $\{1,\ldots,T\}$&Set of decision stages\\
    $\mathcal{N}$ &Set of nodes in the multi-distribution scenario tree\\
    $\mathcal{N}_t\subseteq{\mc{N}}$ & Set of nodes at stage $t$\\
    $\mc{D}_n$ & Set of possible distributions applicable at node $n$\\
    $\mathcal{N}_{nd}\subseteq{\mc{N}}$ & Set of child nodes of node
                                          $n$ if distribution $d\in \mathcal{D}_n$ is enforced\\
    $X_{t}$& Domain of the decision variables at stage $t$\\
    \midrule
    \multicolumn{2}{c}{Parameters}  \\
    \midrule
    $t(n)$ & Stage of node $n$\\
    $a(n)$ & Parent node of node $n$\\
    $\pi_n$ & Probability of node $n$\\
    $r_n\in \mathbb{R}^{N_{t(n)}}$& Coefficients of decision variables $x_n$ in the objective function\\
    $q_{nd}\in \mathbb{R}^1$& Coefficient of decision variable $\delta_{nd}$ in the objective function\\
    $A_n\in \mathbb{R}^{M_{t(n)}\times N_{t(n)}}$& Coefficients of variables $x_n$ in the constraints that connect $x_n$ and $\delta_{nd}$ decisions\\
    $B_{nd}\in \mathbb{R}^{M_{t(n)}\times 1}$&Coefficients of variable $\delta_{nd}$ in the constraints that connect $x_n$ and $\delta_{nd}$ decisions\\
    $C_{a(n)}\in \mathbb{R}^{M_{t(n)}\times N_{t(a(n))}}$& Coefficients of variables $x_{a(n)}$ in the constraints that connect $x_n$ and $\delta_{nd}$ decisions\\
    $D_{a(n),d}\in \mathbb{R}^{M_{t(n)}\times 1}$&Coefficients of variable $\delta_{a(n),d}$ in the constraints that connect $x_n$ and $\delta_{nd}$ decisions\\
    $h_n\in \mathbb{R}^{M_{t(n)}}$& Right-hand-side coefficients of the constraints that connect $x_n$ and $\delta_{nd}$ decisions\\
    $\Theta_n \in \mathbb{R}^1$& Terminal value of the decisions following node $n\in\mc{N}_T$\\
    \midrule
    \multicolumn{2}{c}{Variables}  \\
    \midrule
    $x_n\in \mathbb{R}^{N_{t(n)}}$ & Decisions made at node $n$\\
    $\delta_{nd}\in\{0,1\}$& Decision on whether to apply probability distribution $d$ at node $n$\\
    $\theta_n\in \mathbb{R}^1$& Expected value of the decisions made at the nodes descending from $n$\\
    \bottomrule
  \end{longtable}
  
  \section{A big-$M$ reformulation}\label{sec:app:bigMref}
  In this appendix a big-$M$ reformulation that linearizes model \eqref{eq:eumsp} is introduced.
  In addition to the notation introduced in \Cref{sec:MSPEU}, let $M_{nd}\in \mathbb{R}^1$ be a suitably high constant.
  The linearized EUMSP is thus
  \begin{subequations}\label{eq:ref:eumsp}
    \begin{align}
      \label{eq:ref:eumsp:obj} \max~& r_0^Tx_{0}+\sum_{\mathclap{d\in \mathcal{D}_0}}q_{0d}\delta_{0d}+\theta_0  \\
      \label{eq:ref:eumsp:c1}  \text{s.t.}&\sum_{\mathclap{d\in \mathcal{D}_n}}\delta_{nd} = 1 &n\in \mathcal{N},\\
      \label{eq:ref:eumsp:c2}   &A_nx_n + \sum_{\mathclap{d\in \mathcal{D}_n}}B_{nd}\delta_{nd}+ C_{a(n)}x_{a(n)}   + \sum_{\mathclap{d\in \mathcal{D}_{a(n)}}}D_{a(n),d}\delta_{a(n),d}= h_n& n\in \mathcal{N},\\
      \label{eq:ref:eumsp:c3}   &\theta_n\leq \sum_{\mathclap{m\in \mathcal{N}_{nd}}}\pi_m(r_m^Tx_{m}+\sum_{\mathclap{d\in\mathcal{D}_m}}q_{md}\delta_{md}+\theta_m) + M_{nd} (1 - \delta_{nd})& n\in \mathcal{N}\setminus{\mathcal{N}_T}, d\in \mathcal{D}_n\\
      \label{eq:ref:eumsp:c4}    &\theta_{n} = \Theta_n & n\in \mathcal{N}_{T},\\
      \label{eq:ref:eumsp:c5} & x_n \in X_{t(n)} & n\in \mathcal{N},\\
      \label{eq:ref:eumsp:c6} & \delta_{nd} \in \{0,1\} & n\in \mathcal{N},d\in \mathcal{D}_n,\\
      \label{eq:ref:eumsp:c7} & \theta_{n} \in \mathcal{R} & n\in \mathcal{N}.
    \end{align}
  \end{subequations}
  Note, particularly, that constraints \eqref{eq:ref:eumsp:c3} are equivalent to \cref{eq:eumsp:c3}. Consider a given node $n$, other than a leaf node.
  Observe that only for one distribution $d$ there will be a $\delta_{nd}$ which takes value one at $n$ (see \cref{eq:ref:eumsp:c1}).
  For the same $n$ and for the same $d$, the second term on the right-hand-side of \cref{eq:ref:eumsp:c3} will be zero (i.e., the big-$M$ will not be enforced), and the resulting right-hand-side will be the most binding among the $|\mathcal{D}_n|$ constraints for node $n$. Since we are maximizing, at optimality $\theta_n$ will take value of the expectation according to the distribution $d$ for which $\delta_{nd}=1$, as it happens in model \cref{eq:eumsp}.

  \section{Finding big-$M$ values}\label{sec:app:bigMs}
  An efficient implementation of model \cref{eq:ref:eumsp} requires tight big-$M$ values.
  Observe that, for $n\in\mc{N}\setminus\mc{N}_T$ and $d\in\mc{D}_n$,
  constant $M_{nd}$ must be a valid upper bound for constraints \eqref{eq:ref:eumsp:c3}, that is:
  $$\theta_{n} - \sum_{\mathclap{m\in \mathcal{N}_{nd}}}\pi_m(r_m^Tx_{m}+\sum_{\mathclap{d\in\mathcal{D}_m}}q_{md}\delta_{md}+\theta_m) \leq M_{nd} $$
  Let us introduce $\phi_{nd}$ to represent the expectation at the children of node $n$ for distribution $d$, that is:
  $$\phi_{nd}  = \sum_{\mathclap{m\in \mathcal{N}_{nd}}}\pi_m(r_m^Tx_{m}+\sum_{\mathclap{d\in\mathcal{D}_m}}q_{md}\delta_{md}+\theta_m)$$
  Consider the numerical example shown in \Cref{tab:bigMs:example} for a given node $n$ and three possible distributions $d$. The table reports
  the highest and lowest values the expectation at the following stage can take for each possible distribution, and the corresponding value of $\theta_n$
  should a specific distribution be chosen. When choosing $M_{n,d_1}$, notice that the maximum value $\theta_n$ can reach for other distributions is $9$,
  and that the least value $\phi_{n,d_1}$ can reach is $5$. Therefore, we can set $M_{n,d_1}=9-5$. In fact, if distribution $d_2$ is chosen, $\theta_n$ will be
  at most $9$ and $\phi_{n,d_1}$ at least $5$, thus adding $4$ to $\phi_{n,d_1}$ will ensure that $\theta_n$ is correctly set to $9$. Similarly, we choose $M_{n,d_2}=6$
  as the highest value $\theta_n$ can take if $d_2$ is not selected is $10$, while the least value of $\phi_{n,d_2}$ is $4$. Finally, with a similar reasoning we can set $M_{n,d_3}=7$.

  \begin{table}[h]
    \centering
    \caption{Numerical example for the calculation of constants $M_{nd}$}
    \label{tab:bigMs:example}
    \begin{tabular}{ccccc}
      \toprule
      &&\multicolumn{2}{|c|}{$\phi_{nd}$}&\\
      $d$&$\theta_n$&Max&Min&$M_{nd}$\\
      \midrule
      $d_1$&10&10&5&$9-5=4$\\
      $d_2$&9&9&4&$10-4=6$\\
      $d_3$&8&8&3&$10-3=7$\\
      \bottomrule    
    \end{tabular}
    
  \end{table}
  From the example in \Cref{tab:bigMs:example} we understand that finding values for $M_{nd}$ amounts to finding highest values for $\phi_{nd}$ and differences $\theta_n - \phi_{nd}$.
  In what follows we illustrate how these values can be found for $t=T-1$ in \Cref{sec:bigM:T-1} and $t=T-2,\ldots,1$ in \Cref{sec:bigM:T-2}.

  \subsection{Big-$M$ values for stages $t=T-1$}\label{sec:bigM:T-1}

  We start at the second-last stage, $t=T-1$. Our task is that of finding, for each node $\bar{n}\in\mathcal{N}_{T-1}$ and for each distribution $\bar{d}\in \mathcal{D}_{\bar{n}}$,
  a constant $M_{\bar{n}\bar{d}}$ which is slightly higher than the highest difference $\theta_{\bar{n}}-\phi_{\bar{n}\bar{d}}$, where again
  $$\phi_{\bar{n}\bar{d}} = \sum_{\mathclap{m\in \mathcal{N}_{\bar{n}\bar{d}}}}\pi_m(r_m^Tx_{m}+\sum_{\mathclap{d\in\mathcal{D}_m}}q_{md}\delta_{md}+\theta_m)$$
  Now, the highest difference can be found solving the following optimization problem:
  \begin{subequations}\label{eq:bigM:T-1}
    \begin{align}
      \label{eq:bigM:T-1:obj} M_{\bar{n}\bar{d}}^*=\max~& \theta_{\bar{n}} - \sum_{m\in \mathcal{N}_{\bar{n}\bar{d}}}\pi_m(r_m^Tx_{m}+\sum_{\mathclap{k\in\mathcal{D}_m}}q_{mk}\delta_{mk}+\Theta_m)  \\
      \label{eq:bigM:T-1:c1}  \text{s.t.}&\sum_{\mathclap{d\in \mathcal{D}_n}}\delta_{nd} = 1 &n\in \mathcal{N},\\
      \label{eq:bigM:T-1:c2}   &A_nx_n + \sum_{\mathclap{d\in \mathcal{D}_n}}B_{nd}\delta_{nd}+ C_{a(n)}x_{a(n)}   + \sum_{\mathclap{d\in \mathcal{D}_{a(n)}}}D_{a(n),d}\delta_{a(n),d}= h_n& n\in \mathcal{N},\\
      \label{eq:bigM:T-1:c4} &\theta_{\bar{n}} \leq \Theta_{\bar{n}\bar{d}}^*,&\\
      \label{eq:bigM:T-1:c5} & x_n \in X_{t(n)} & n\in \mathcal{N},\\
      \label{eq:bigM:T-1:c6} & \delta_{nd} \in \{0,1\} & n\in \mathcal{N},d\in \mathcal{D}_n,\\
      \label{eq:bigM:T-1:c6} & \delta_{\bar{n}\bar{d}} = 0 &
    \end{align}
  \end{subequations}
  Problem \cref{eq:bigM:T-1} consists of finding the feasible solution to problem \eqref{eq:ref:eumsp} which yields the highest value for the left-hand-side of constraint \eqref{eq:ref:eumsp:c3} for $\bar{n}$ and $\bar{d}$.
  The following two elements must be noted in \eqref{eq:bigM:T-1}. Constraints \eqref{eq:ref:eumsp:c3} of the original problem, which determine the correct expectations at the stages before $T-1$, are not included as they are irrelevant for
  stage $T-1$. The second element to note is constraint \eqref{eq:bigM:T-1:c4} which sets an upper bound $\theta_{\bar{n}}$. This upper bound represents the highest value $\theta_{\bar{n}}$ can take for the distributions other than $\bar{d}$.
  This value can, in turn, be obtained solving optimization problems. The highest expectation for stage $T$, given distribution $d'\in{D}_{\bar{n}}$ is the optimal value to problem \eqref{eq:bigM:T-1:maxTheta}:
  \begin{subequations}\label{eq:bigM:T-1:maxTheta}
    \begin{align}
      \label{eq:bigM:T-1:maxTheta:obj} \Phi_{\bar{n}d'}^*=\max~& \sum_{m\in \mathcal{N}_{\bar{n}d'}}\pi_m(r_m^Tx_{m}+\sum_{\mathclap{k\in\mathcal{D}_m}}q_{mk}\delta_{mk}+\Theta_m)  \\
      \label{eq:bigM:T-1:maxTheta:c1}  \text{s.t.}&\sum_{\mathclap{d\in \mathcal{D}_n}}\delta_{nd} = 1 &n\in \mathcal{N},\\
      \label{eq:bigM:T-1:maxTheta:c2}   &A_nx_n + \sum_{\mathclap{d\in \mathcal{D}_n}}B_{nd}\delta_{nd}+ C_{a(n)}x_{a(n)}   + \sum_{\mathclap{d\in \mathcal{D}_{a(n)}}}D_{a(n),d}\delta_{a(n),d}= h_n& n\in \mathcal{N},\\
      \label{eq:bigM:T-1:maxTheta:c3} & x_n \in X_{t(n)} & n\in \mathcal{N},\\
      \label{eq:bigM:T-1:maxTheta:c4} & \delta_{nd} \in \{0,1\} & n\in \mathcal{N},d\in \mathcal{D}_n,\\
      \label{eq:bigM:T-1:maxTheta:c5} &\delta_{\bar{n}d'} = 1 &      
    \end{align}
  \end{subequations}
  Therefore, when calculating $M_{\bar{n}\bar{d}}$, the upper bound $\Theta_{\bar{n}}^*$ in \cref{eq:bigM:T-1:c4} is given by:
  $$\Theta_{\bar{n}\bar{d}}^* = \max_{ d\in \mc{D}_{\bar{n}} : d\neq \bar{d}}\Phi_{\bar{n}d}^*$$

  Clearly, solving problems \eqref{eq:bigM:T-1} and \eqref{eq:bigM:T-1:maxTheta} amounts to solving integer programs of size comparable with the original problem \eqref{eq:ref:eumsp}.
  However, the tightest $\Theta_{\bar{n}\bar{d}}^*$ and $M_{\bar{n}\bar{d}}$ are not necessary, and higher values would still provide correct results.
  A suitable value for $M_{\bar{n}\bar{d}}$ can be obtained by solving any relaxation of problem \eqref{eq:bigM:T-1}, yielding $M_{\bar{n}\bar{d}}^{R}$,
  and \eqref{eq:bigM:T-1:maxTheta}, yielding $\Phi_{\bar{n}\bar{d}}^{R}$ and in turn $\Theta_{\bar{n}\bar{d}}^{R}$. As an example, one might solve the linear programming relaxation of problems \eqref{eq:bigM:T-1}
  and \eqref{eq:bigM:T-1:maxTheta} or, if the size of the problems is excessively high, one might choose to relax constraints \cref{eq:bigM:T-1:c2} and \cref{eq:bigM:T-1:maxTheta:c2} for some stages.
  Finally, since the procedure outlined might return negative values for some $M_{nd}$, we set $M_{nd}=\max\{0,M_{nd}^{R}\}$ to reduce, when possible, high big-$M$ absolute values.
  The procedure is summarized in \Cref{alg:bigM:T-1}.

  \begin{algorithm}[h]
    \caption{Algorithm for calculating $M_{nd}$ for $T-1$}
    \label{alg:bigM:T-1}
    \begin{algorithmic}[1]
      \STATE Input: $\mc{N}$, $\mc{D}_n$ for $n\in \mc{N}$, $\Theta_n$ for $n \in \mc{N}_T$
      \FOR{Node $\bar{n}\in \mc{N}_{T-1}$}
      \FOR{Distribution $d\in \mc{D}_{\bar{n}}$}
      \STATE Calculate $\Phi^{LP}_{\bar{n}d}$ by solving the LP relaxation to problem \cref{eq:bigM:T-1:maxTheta}
      \ENDFOR
      \FOR{Distribution $\bar{d}\in \mc{D}_{\bar{n}}$}
      \STATE In constraint \cref{eq:bigM:T-1:c4} set $\Theta_{\bar{n}\bar{d}}^*=\max_{d\in \mc{D}_{\bar{n}}:d\neq \bar{d}}\Phi_{\bar{n}d}^{LP}$
      \STATE Calculate $M_{\bar{n}\bar{d}}^R$ by solving a suitable relaxation to problem \cref{eq:bigM:T-1}
      \STATE Set $M_{\bar{n}\bar{d}} = \max\{0,M_{\bar{n}\bar{d}}^R\}$
      \ENDFOR
      \ENDFOR
      \RETURN $M_{nd}$ for $n \in \mc{N}_{T-1}$ and $d\in \mc{D}_n$.
    \end{algorithmic}
  \end{algorithm}

  \subsection{Big-$M$ values for stages $t=T-2,\ldots,1$}\label{sec:bigM:T-2}
  Once constants $M_{nd}$ are available for every $n\in \mc{N}_{T-1}$ and $d\in\mc{D}_n$ we can proceed in a similar way to calculate big-$M$s for stages $T-2,\ldots,1$.
  Given a stage $\bar{t}\in \{T-2,\ldots,1\}$, a node at that stage, $\bar{n}\in \mc{N}_{\bar{t}}$, and distribution available at that node $\bar{d}\in \mc{D}_{\bar{n}}$,
  the tightest value of constant $M_{\bar{n}\bar{d}}$, namely $M_{\bar{n}\bar{d}}^*$, is

  \begin{footnotesize}
    \begin{subequations}\label{eq:bigM:T-2}
      \begin{align}
        \label{eq:bigM:T-2:obj} M_{\bar{n}\bar{d}}^*=\max~& \theta_{\bar{n}} - \sum_{m\in \mathcal{N}_{\bar{n}\bar{d}}}\pi_m(r_m^Tx_{m}+\sum_{\mathclap{k\in\mathcal{D}_m}}q_{mk}\delta_{mk}+\theta_m)  \\
        \label{eq:bigM:T-2:c1}  \text{s.t.}&\sum_{\mathclap{d\in \mathcal{D}_n}}\delta_{nd} = 1 &n\in \mathcal{N},\\
        \label{eq:bigM:T-2:c2}   &A_nx_n + \sum_{\mathclap{d\in \mathcal{D}_n}}B_{nd}\delta_{nd}+ C_{a(n)}x_{a(n)}   + \sum_{\mathclap{d\in \mathcal{D}_{a(n)}}}D_{a(n),d}\delta_{a(n),d}= h_n& n\in \mathcal{N},\\
        \label{eq:bigM:T-2:c3} &\theta_{n} \geq \sum_{m\in \mathcal{N}_{nd}}\pi_m(r_m^Tx_{m}+\sum_{\mathclap{k\in\mathcal{D}_m}}q_{mk}\delta_{mk}+\theta_m) - M_{nd} (1 - \delta_{nd})& t=\bar{t}+1,\ldots,T-1,n\in \mc{N}_t,d\in \mathcal{D}_{n},\\
        \label{eq:bigM:T-2:c4} &\theta_{n} = \Theta_{n}&n\in \mc{N}_T\\
        \label{eq:bigM:T-2:c5} &\theta_{\bar{n}} \leq \Theta_{\bar{n}\bar{d}}^*,&\\
        \label{eq:bigM:T-2:c6} & x_n \in X_{t(n)} & n\in \mathcal{N},\\
        \label{eq:bigM:T-2:c7} & \delta_{nd} \in \{0,1\} & n\in \mathcal{N},d\in \mathcal{D}_n,\\
        \label{eq:bigM:T-2:c8} &\delta_{\bar{n}\bar{d}} = 0 &
      \end{align}
    \end{subequations}
  \end{footnotesize}

  Problem \cref{eq:bigM:T-2} consists of finding the feasible solution to problem \cref{eq:ref:eumsp} which yields the highest value of the left-hand-side of constraint \cref{eq:ref:eumsp:c3} for node $\bar{n}$ and distribution $\bar{d}$.
  Notice that, unlike in problem \cref{eq:bigM:T-1}, problem \cref{eq:bigM:T-2} includes constraints \cref{eq:bigM:T-2:c3} which are necessary to ensure that $\theta_n$ values are set to the lowest expectation for all stages between $\bar{t}$ and $T-1$. Note that constant $M_{nd}$ in constraints \cref{eq:bigM:T-2:c3} is also an upper bound to the quantity $\theta_{n} - \sum_{m\in \mathcal{N}_{nd}}\pi_m(r_m^Tx_{m}+\sum_{k\in\mathcal{D}_m}q_{mk}\delta_{mk}+\theta_m)$ and can be thus set to the quantities determined at previous iterations. Also in this case, $\Theta_{\bar{n}}^*$ in \cref{eq:bigM:T-2:c5} represents the highest possible value $\theta_{\bar{n}}$
  can take for distributions other than $\bar{d}$. The highest expectation for stage $\bar{t}+1$, given distribution $d'\in{D}_{\bar{n}}$ is the optimal value to problem \eqref{eq:bigM:T-2:maxTheta}
  \begin{footnotesize}
    \begin{subequations}\label{eq:bigM:T-2:maxTheta}
      \begin{align}
        \label{eq:bigM:T-2:maxTheta:obj} \Phi_{\bar{n}d'}^*=\max~& \sum_{m\in \mathcal{N}_{\bar{n}d'}}\pi_m(r_m^Tx_{m}+\sum_{\mathclap{k\in\mathcal{D}_m}}q_{mk}\delta_{mk}+\theta_m)  \\
        \label{eq:bigM:T-2:maxTheta:c1}  \text{s.t.}&\sum_{\mathclap{d\in \mathcal{D}_n}}\delta_{nd} = 1 &n\in \mathcal{N},\\
        \label{eq:bigM:T-2:maxTheta:c2}   &A_nx_n + \sum_{\mathclap{d\in \mathcal{D}_n}}B_{nd}\delta_{nd}+ C_{a(n)}x_{a(n)}   + \sum_{\mathclap{d\in \mathcal{D}_{a(n)}}}D_{a(n),d}\delta_{a(n),d}= h_n& n\in \mathcal{N},\\
        \label{eq:bigM:maxTheta:T-2:c3} &\theta_{n} \leq \sum_{m\in \mathcal{N}_{nd}}\pi_m(r_m^Tx_{m}+\sum_{\mathclap{k\in\mathcal{D}_m}}q_{mk}\delta_{mk}+\theta_m)+ M_{nd} (1 - \delta_{nd})& t=\bar{t}+1,\ldots,T-1,n\in \mc{N}_t,d\in \mathcal{D}_{n},\\
        \label{eq:bigM:maxTheta:T-2:c4} &\theta_{n} = \Theta_{n}&n\in \mc{N}_T\\
        \label{eq:bigM:T-2:maxTheta:c5} & x_n \in X_{t(n)} & n\in \mathcal{N},\\
        \label{eq:bigM:T-2:maxTheta:c6} & \delta_{nd} \in \{0,1\} & n\in \mathcal{N},d\in \mathcal{D}_n,\\
        \label{eq:bigM:T-2:maxTheta:c7} &\delta_{\bar{n}d'} = 1 &      
      \end{align}
    \end{subequations}
  \end{footnotesize}
  Therefore, the upper bound $\Theta_{\bar{n}\bar{d}}^*$ in \cref{eq:bigM:T-2:c5} is given by:
  $$\Theta_{\bar{n}}^* = \max_{d\in \mc{D}_{\bar{n}}:d\neq \bar{d}}\Phi_{\bar{n}d}^*$$
  Similarly to \Cref{sec:bigM:T-1}, calculating the optimal $M_{\bar{n}\bar{d}}^*$ is cumbersome as well as not strictly necessary. Therefore, any computationally suitable relaxation to problems \cref{eq:bigM:T-2}
  and \cref{eq:bigM:T-2:maxTheta} can be adopted. The procedure for obtaining constants $M_{nd}$ for stages $T-2,\ldots,1$ is sketched in \Cref{alg:bigM:T-2}.

  \begin{algorithm}[h]
    \caption{Algorithm for calculating $M_{nd}$ for $\bar{t} = T-2,\ldots,1$}
    \label{alg:bigM:T-2}
    \begin{algorithmic}[1]
      \STATE Input: $\mc{N}$, $\mc{D}_n$ for $n\in \mc{N}$, $\Theta_n$ for $n \in \mc{N}_T$, $M_{nd}$ for $n\in\mc{N}_t$, $t=\bar{t}+1,\ldots,T-1$, $d\in {D}_n$.
      \FOR{Node $\bar{n}\in \mc{N}_{\bar{t}}$}
      \FOR{Distribution $d\in \mc{D}_{\bar{n}}$}
      \STATE Calculate $\Phi^{LP}_{\bar{n}d}$ by solving the LP relaxation to problem \cref{eq:bigM:T-2:maxTheta}
      \ENDFOR
      \FOR{Distribution $\bar{d}\in \mc{D}_{\bar{n}}$}
      \STATE In constraint \cref{eq:bigM:T-2:c5} set $\Theta_{\bar{n}\bar{d}}^*=\max_{d\in \mc{D}_{\bar{n}}:d\neq \bar{d}}\Phi_{\bar{n}d}^{LP}$
      \STATE Calculate $M_{\bar{n}\bar{d}}^R$ by solving a relaxation to problem \cref{eq:bigM:T-2}
      \STATE Set $M_{\bar{n}\bar{d}} = \max\{0,M_{\bar{n}\bar{d}}^R\}$
      \ENDFOR
      \ENDFOR
      \RETURN $M_{nd}$ for $n \in \mc{N}_{\bar{t}}$ and $d\in \mc{D}_n$.
    \end{algorithmic}
  \end{algorithm}

  \section{Computation time for finding big-$M$ values for the FTCP }\label{app:bigMtime}
  
  \begin{longtable}{lc|ccc}
    \caption{Average elapsed time in seconds for the computations of big-$M$ values using the procedure in \Cref{sec:cs:bigM}. $S$ indicates the number of realizations describing each distribution.}\label{tab:bigMs}\\
  \toprule
    Team &    $S$ & $|\mc{I}|=3$ & $|\mc{I}|=4$& $|\mc{I}|=5$\\
    \midrule
     ARS &         4 &     11.216 &     36.161 &   66.495 \\
 ASV &         4 &     13.387 &     37.614 &   84.259 \\
 CAR &         4 &     17.873 &     44.370 &   96.844 \\
 CHE &         4 &     16.862 &     52.053 &  149.683 \\
 CRP &         4 &     17.100 &     47.362 &  129.843 \\
 EVE &         4 &     10.921 &     36.309 &   93.248 \\
 FUL &         4 &     12.467 &     39.589 &   67.543 \\
 HUL &         4 &     11.928 &     40.945 &   86.544 \\
 LIV &         4 &     11.384 &     32.654 &   84.933 \\
 MAC &         4 &     12.999 &     32.460 &   71.190 \\
 MAU &         4 &     14.093 &     38.464 &   95.931 \\
 NEC &         4 &     11.079 &     32.661 &   87.461 \\
 NOR &         4 &      9.897 &     28.554 &   66.370 \\
 SOU &         4 &     10.861 &     23.151 &   66.977 \\
 STO &         4 &     11.009 &     37.650 &   68.867 \\
 SUN &         4 &     18.335 &     36.425 &   85.658 \\
 SWA &         4 &     12.770 &     32.966 &   99.143 \\
 TOT &         4 &     11.475 &     28.322 &   69.286 \\
 WBA &         4 &     10.644 &     28.465 &   85.157 \\
 WHU &         4 &     17.390 &     45.385 &   78.934 \\
    ARS &         5 &     23.028 &     59.372 &  165.883 \\
 ASV &         5 &     21.531 &     62.066 &  167.843 \\
 CAR &         5 &     24.472 &     74.075 &  241.720 \\
 CHE &         5 &     28.144 &     84.103 &  301.168 \\
 CRP &         5 &     21.557 &     74.754 &  257.036 \\
 EVE &         5 &     21.673 &     72.323 &  183.927 \\
 FUL &         5 &     19.659 &     50.441 &  142.312 \\
 HUL &         5 &     23.749 &     51.421 &  176.955 \\
 LIV &         5 &     18.195 &     40.295 &  150.142 \\
 MAC &         5 &     20.799 &     47.659 &  117.553 \\
 MAU &         5 &     31.487 &     64.327 &  155.962 \\
 NEC &         5 &     35.180 &     64.231 &  175.313 \\
 NOR &         5 &     26.893 &     51.132 &  146.938 \\
 SOU &         5 &     18.282 &     43.805 &  120.134 \\
 STO &         5 &     21.236 &     54.687 &  151.571 \\
 SUN &         5 &     29.861 &     66.360 &  187.473 \\
 SWA &         5 &     30.952 &     59.528 &  141.202 \\
 TOT &         5 &     28.963 &     56.107 &  128.371 \\
 WBA &         5 &     23.609 &     54.868 &  139.128 \\
    WHU &         5 &     27.754 &     63.924 &  176.582 \\
    \midrule
   Avg &     &     19.018 &     48.176 &  129.039 \\
    \bottomrule
  \end{longtable}
  

\end{appendices}

\end{document}